\documentclass{article}

\usepackage{
 amsmath,
 amsxtra,
 amsthm,
 amssymb,
 etex,
 mathrsfs,
 }
\usepackage[all]{xy}

\newtheorem{theorem}{Theorem}[section]
\newtheorem{lemma}[theorem]{Lemma}
\newtheorem{conjecture}[theorem]{Conjecture}
\newtheorem{proposition}[theorem]{Proposition}
\newtheorem{corollary}[theorem]{Corollary}

\theoremstyle{definition}
\newtheorem{defn}[theorem]{Definition}
\newtheorem{remark}[theorem]{Remark}

\newcommand{\bd}{\begin{defn}}
\newcommand{\ed}{\end{defn}}
\newcommand{\bl}{\begin{lemma}}
\newcommand{\el}{\end{lemma}}
\newcommand{\bp}{\begin{proposition}}
\newcommand{\ep}{\end{proposition}}
\newcommand{\bt}{\begin{theorem}}
\newcommand{\et}{\end{theorem}}
\newcommand{\bc}{\begin{corollary}}
\newcommand{\ec}{\end{corollary}}
\newcommand{\br}{\begin{remark}}
\newcommand{\er}{\end{remark}}
\newcommand{\ba}{\begin{array}}
\newcommand{\ea}{\end{array}}
\newcommand{\bpf}{\begin{proof}}
\newcommand{\epf}{\end{proof}}

\newcommand{\Z}{\mathbb{Z}}
\newcommand{\Q}{\mathbb{Q}}
\newcommand{\Zp}{\mathbb{Z}_{p}}
\newcommand{\Qp}{\mathbb{Q}_{p}}
\newcommand{\Op}{\mathcal{O}}
\newcommand{\Ep}{E_{p^{\infty}}}
\newcommand{\al}{\alpha}
\newcommand{\be}{\beta}
\newcommand{\Ga}{\Gamma}
\newcommand{\ga}{\gamma}

\newcommand{\la}{\lambda}

\newcommand{\Si}{\Sigma}

\newcommand{\mG}{\mathcal{G}}
\newcommand{\mH}{\mathcal{H}}

\DeclareMathOperator{\Sel}{Sel} \DeclareMathOperator{\Gal}{Gal}
\DeclareMathOperator{\Hom}{Hom} \DeclareMathOperator{\rank}{rank}
\DeclareMathOperator{\corank}{corank}
\DeclareMathOperator{\Ext}{Ext} 

\newcommand{\reg}{\mathrm{reg}}
\newcommand{\res}{\mathrm{res}}

\newcommand{\ord}{\mathrm{ord}}
\newcommand{\cyc}{\mathrm{cyc}}
\newcommand{\ch}{\mathrm{char}}

\newcommand{\tw}{\mathrm{tw}}
\newcommand{\Ak}{\mathrm{Ak}}
\newcommand{\M}{\mathfrak{M}}

\newcommand{\mK}{\mathcal{K}}
\newcommand{\mL}{\mathcal{L}}
\newcommand{\mF}{\mathcal{F}}

\newcommand{\ot}{\otimes}
\newcommand{\ilim}{\displaystyle \mathop{\varinjlim}\limits}
\newcommand{\plim}{\displaystyle \mathop{\varprojlim}\limits}
\newcommand{\im}{\mathrm{im}\,}
\newcommand{\coker}{\mathrm{coker}\,}

\newcommand{\lra}{\longrightarrow}
\newcommand{\tha}{\twoheadrightarrow}

\newcommand{\ps}[1]{[[ #1 ]]}

  \DeclareFontFamily{U}{wncy}{}
  \DeclareFontShape{U}{wncy}{m}{n}{<->wncyr10}{}
  \DeclareSymbolFont{mcy}{U}{wncy}{m}{n}
  \DeclareMathSymbol{\sha}{\mathord}{mcy}{"58}

\numberwithin{equation}{section}

\begin{document}
\title{On order of vanishing of characteristic elements}
 \author{
  Meng Fai Lim\footnote{School of Mathematics and Statistics $\&$ Hubei Key Laboratory of Mathematical Sciences,
Central China Normal University, Wuhan, 430079, P.R.China.
 E-mail: \texttt{limmf@ccnu.edu.cn}} }
\date{}
\maketitle

\begin{abstract} \footnotesize
\noindent
Let $p$ be a fixed odd prime.
Let $E$ be an elliptic curve defined over a number field $F$ with either good ordinary reduction or multiplicative reduction at each prime of $F$ above $p$.  We shall study the characteristic element of the Selmer group of $E$ over a $p$-adic Lie extension. In particular, we relate the order of vanishing of these characteristic element evaluated at Artin representations to the Selmer coranks and their twists in the intermediate subextensions of the $p$-adic Lie extension.

\medskip
\noindent Keywords and Phrases: Characteristic element, Selmer groups, order of vanishing.

\smallskip
\noindent Mathematics Subject Classification 2020: 11G05, 11R23, 11S25.
\end{abstract}

\section{Introduction}

Over the last two decades or so, Selmer groups over non-commutative $p$-adic Lie extensions have been much studied by many. Initial attempts were
for proving some form of structure theorem for finitely generated (torsion) modules over Iwasawa algebra of compact $p$-adic Lie groups in hope of attaching characteristic elements to these Selmer groups.
As this turned out to be too difficult (see \cite{CSSAlg} and especially the introduction in \cite{CF+}),
Venjakob \cite{V05} came up with an idea of bypassing structure theorems using algebraic $K$-theory of Iwasawa algebras. Building on Venjakob's innovation, Coates-Fukaya-Kato-Sujatha-Venjakob \cite{CF+, FK} went on to formulate a non-commutative analogue of the Iwasawa main conjecture (also see \cite{Bu15, BV, DL, DD}). More precisely, they were able to attach characteristic elements to Selmer groups over a $p$-adic Lie extension (under certain hypothesis). These characteristic elements live in a localized $K_1$-group by virtue of their definition, and are conjectured to interpolate the special values of the complex $L$-functions of the elliptic curve twisted by Artin representations of the Galois group of the said $p$-adic Lie extension. In this paper, we shall examine these characteristic elements and their relation with the Selmer coranks in the intermediate subextensions of the $p$-adic Lie extension. For the remainder of the introductory section, we shall say a bit more, leaving details to the body of the paper.

Throughout, $p$ will always denote an odd prime. For simplicity, we assume in this introduction that our elliptic curve $E$ is defined over $\Q$. The elliptic curve is further assumed to have either good ordinary reduction or multiplicative reduction (possibly split or non-split) at the prime $p$. Write $\Q^\cyc$ for the cyclotomic $\Zp$-extension of $\Q$. Let $F_\infty$ be a Galois extension of $\Q$ containing $\Q^\cyc$ with Galois group $\Gal(F_{\infty}/\Q)$ being a compact $p$-adic Lie group with no $p$-torsion.
Write $G=\Gal(F_\infty/\Q)$, $H=\Gal(F_\infty/\Q^\cyc)$ and $\Ga = G/H = \Gal(\Q^\cyc/\Q)$. Denote by $\M_H(G)$ the category of finitely generated $\Zp\ps{G}$-module $M$ with the property that $M/M[p^\infty]$ is finitely generated over $\Zp\ps{H}$. Under the assumption that the dual Selmer group $X(E/F_\infty)$ of $E$ over the $p$-adic Lie extension $F_\infty$ lies in $\M_H(G)$ (see \cite{CF+,CS12}; also see Conjecture \ref{MHG conj}), one can attach a characteristic element $\xi_E$ to $X(E/F_\infty)$ in the sense of Coates \textit{et al} (see \cite{CF+} or body of our paper). Let $\rho$ be an Artin representation of $\Gal(F_\infty/\Q)$ with coefficients in $\Op$, where $\Op$ is the ring of integers of a finite extension of $\Qp$. In their paper, Coates \textit{et al} introduced a notion of evaluating $\xi_E$ at $\rho$ which yields an element in the field of fraction of the ring $\Op\ps{\Ga}$. This latter element will be denoted by $\Phi_{\rho}(\xi_E)$. Under a fixed choice of topological generator of $\Ga$, we may view $\Phi_{\rho}(\xi_E)$ as an element in the field of fraction of $\Op\ps{T}$, where $\Op\ps{T}$ is the power series ring in one variable. Therefore, it makes sense to speak of
\[ \ord_{T=0}\big(\Phi_{\rho}(\xi_E)\big) \]
(cf. \cite{Bu15}; also see Definition \ref{BurnsDef}).

Now, if $L$ is a finite extension of $\Q$ contained in $F_\infty$, write $\reg_L$ for the regular representation of $\Gal(L/\Q)$. Our first main result is concerned with the order of vanishing of the characteristic elements evaluated at these regular representations (see Theorem \ref{main: pro-p} below for a more precise and general version of the theorem).

\bt
Retain notation as above. Assume that $X(E/F_\infty)\in \M_H(G)$. Suppose that $L$ is a finite Galois extension of $\Q$ contained in $F_\infty$ which satisfies the following three statements.
\begin{enumerate}
\item[$(a)$] $F_\infty/L$ is a pro-$p$ extension.
\item[$(b)$] $X(E/L^\cyc)$ satisfies a semisimple conjecture of Greenberg $($see Conjecture \ref{semisimple conj}$)$.
\item[$(c)$] $L$ satisfies property $\mathbf{(M_p)}$ $($see Definition \ref{splitmultdef}$)$.
\end{enumerate}
If $\xi_E$ is a characteristic element of $X(E/F_\infty)$, we have
\[\displaystyle\ord_{T=0} ~\big(\Phi_{\reg_L}(\xi_E)\big) = \rank_{\Zp}\big(X(E/L)\big) + m_p(L),\]
where $m_p(L)$ is an integer defined as in Definition \ref{splitmultdef}.
\et

We mention that the term $m_p(L)$ only comes into play when $E$ has split multiplicative reduction, and there do exist situations, where $m_p(L)$ can be non-zero which was first observed in \cite{Ze11}. However, we also note that even in the presence of split multiplicative reduction prime, this quantity can still be zero. For instance, if $F_\infty = \Q(\mu_{p^\infty}, \sqrt[p^\infty]{m})$ for some $p$-power-free integer $m>1$, then $m_p(L)$ is always zero (see Lemma \ref{akashi zero order split mult 2}).

It is natural to ask if one can obtain a similar result for Artin representations which are not regular representations. This is our next result. To explain, we need to introduce more notation. Suppose that the Artin representation $\rho$ is irreducible. We then write $\mF$ for any finite Galois extension of $\Q$ such that $\rho$ factors through $\Gal(\mF/\Q)$. Let $s_{E,\rho}$ denote the number of copies of $W_\rho\ot_{\Op}\bar{\Q}_p$ occurring in $X(E/\mF)\ot_{\Zp}\bar{\Q}_p$. We shall also write $X(\tw_{\rho}(E)/\Q^\cyc)$ for the twisted Selmer group (see Subsection \ref{twisted Selmer Subsec} for the definition). Our second result is then as follows (again, see Theorem \ref{main: Artin twist} for a more general version).

\bt \label{main: Artin twist intro}
Let $E$ be an elliptic curve defined over $\Q$ with good ordinary reduction at $p$. Let $F_\infty$ be a $p$-adic Lie extension of $\Q$ with $\Gal(F_\infty/\Q)$ having no $p$-torsion. Write $H=\Gal(F_\infty/\Q)$. Suppose that all of the following statements are valid.
\begin{enumerate}
\item[$(a)$] $X(E/F_\infty)\in \M_H(G)$.
\item[$(b)$] $X(\tw_{\rho}(E)/\Q^\cyc)$ satisfies a semisimple conjecture in the sense of Conjecture \ref{semisimple conj2}.
\item[$(c)$] For every open subgroup $H'$ of $H$, $H^i(H', \Ep(F_\infty))$ is finite for all $i\geq 1$.
 \item[$(d)$] For every prime $w$ of $F^\cyc$ dividing $p$, and each open subgroup $H_w'$ of $H_w$, $H^i(H'_w, \widetilde{E}_{p^\infty}(F_\infty))$ is finite for all $i\geq 1$.
\end{enumerate}
Let $\xi_E$ be a characteristic element of $X(E/F_\infty)$ and $\rho$ an irreducible Artin representation of $\Gal(F_\infty/\Q)$. Then
we have
\[\displaystyle\ord_{T=0} ~\big(\Phi_{\rho}(\xi_E)\big) = s_{E,\rho}.\]
\et

Note that the above theorem applies for a trivial representation. Thus, this removes condition (a) in Theorem 1.1 at the expenses of the extra hypotheses (c) and (d) in Theorem 1.2. Thankfully, these latter hypotheses are known to be satisfied in many $p$-adic Lie extensions (see Remark \ref{hypotheses cd}) and so they are rather mild.

Readers would have noticed that Theorem 1.2 requires $E$ to have good ordinary reduction at primes above $p$. The reason behind this is due to that this theorem relies on certain local cohomology computation of the Artin twist of $\widetilde{E}_{p^\infty}(F_\infty)$ (see Lemma \ref{akashi ord p twist}). At present, we are not able to perform such an analogous computation for an elliptic curve with multiplicative reduction at some prime above $p$. Nevertheless, it would be of interest to have a variant of Theorem 1.2 in the multiplicative reduction context. For this, we are able to obtain such a variant over a False-Tate extension by relying on the (nice enough) representation theory of the intermediate Galois groups of the False-Tate extension. To at least describe this result (again, we refer readers to the body of the paper for the detailed verification), we let $m$ be a $p$-powerfree integer. Set $F_\infty = \Q(\mu_{p^\infty}, \sqrt[p^\infty]{m})$. For $n\geq 1$, let $\rho_n$ denote the representation of $G$ obtained by inducing
a character of exact order $p^n$ of $\Gal\big(\Q(\mu_{p^n}, \sqrt[p^n]{m})/\Q(\mu_{p^n})\big)$ to $\Gal\big(\Q(\mu_{p^n}, \sqrt[p^n]{m})/\Q\big)$. The above asserted result is as follow.

\bt $($see Theorem \ref{ArtinFT}$)$
Let $E$ be an elliptic curve defined over a number field $\Q$ which has either good ordinary reduction or multiplicative reduction at $p$. Suppose that $X(E/F_\infty)\in \M_H(G)$. Assume that the semisimple conjecture $($Conjecture \ref{semisimple conj}$)$ is valid for $X(E/L^\cyc)$, where $L = \Q(\mu_{p^n}, \sqrt[p^n]{m})$.

If $\xi_E$ is a characteristic element of $X(E/F_\infty)$, we have
\[\displaystyle\ord_{T=0} ~\big(\Phi_{\rho_n}(\xi_E)\big) = s_{E,\rho_n}.\]
\et

We end the section giving an outline of the paper. In Section \ref{Alg Sec}, we collect several results on Iwasawa algebras and their modules. In particular, we introduce the procedure of attaching characteristic elements to a certain class of Iwasawa modules which the Selmer groups are expected to belong to. In Section \ref{Elliptic local field}, we review certain properties of elliptic curves over a local field to prepare for the subsequent discussion in the paper. This is followed by the introduction of Selmer groups in Section \ref{Sel Sec}. We also collect several properties of these Selmer groups over a cyclotomic $\Zp$-extension and over a $p$-adic Lie extension. In particular, we introduce the semisimple conjecture of Greenberg (Conjecture \ref{semisimple conj}) and supply a simple criterion that we know of in verifying such a conjecture (see Lemma \ref{semisimple lemma}). We should mention that this criterion is inspired by the work of Wuthrich \cite{Wu}. Section \ref{main results} is where we establish our first result (see Theorem \ref{main: pro-p}). The next two results on the order of vanishing at Artin representations are discussed and proved in Section \ref{twisted Section}. In Section \ref{examples section}, we discuss some classes of examples, where we can calculate the order of vanishing unconditionally. Building on these calculations, we show that the order of vanishing of the characteristic element is bounded above by the order of zero of the Hasse-Weil $L$-function at $s=1$ for these classes of elliptic curves (see Propositions \ref{CS DD}, \ref{Lee2} and \ref{LeeNonsplit}). In some cases, the inequality can even be shown to be an equality (see Corollary \ref{DarmonTian}).

\subsection*{Acknowledgement}
We thank the anonymous referee for the many helpful comments and suggestions.
This research is supported by the
National Natural Science Foundation of China under Grant No. 11771164 and the Fundamental Research Funds for the Central Universities of CCNU under grant CCNU20TD002.

\section{Algebraic preliminaries} \label{Alg Sec}

\subsection{Iwasawa algebras}

Throughout, $G$ will always denote a compact $p$-adic Lie group without $p$-torsion. The Iwasawa algebra (or completed group algebra) of $G$ over $\Zp$ is defined by
 \[ \Zp\ps{G} = \plim_U \Zp[G/U], \]
where $U$ runs over the open normal subgroups of $G$ and the inverse
limit is taken with respect to the canonical projection maps. It is well-known that $\Zp\ps{G}$ is
a Noetherian Auslander regular ring (cf.\ \cite[Theorem 3.26]{V02}).

 In the event that $G$ is pro-$p$, the ring $\Zp\ps{G}$ is local and has no zero divisors (see \cite{Neu}).
Hence it admits a skew field $Q(G)$ which is flat
over $\Zp\ps{G}$ (see \cite[Chapters 6 and 10]{GW} or \cite[Chapter
4, \S 9 and \S 10]{Lam}). Thanks to this property, one can define the $\Zp\ps{G}$-rank of a finitely generated $\Zp\ps{G}$-module $M$ by setting
\[\rank_{\Zp\ps{G}}(M)  = \dim_{Q(G)} (Q(G)\ot_{\Zp\ps{G}}M).\]
The $\Zp\ps{G}$-module $M$ is then said to be torsion  if $\rank_{\Zp\ps{G}} (M) = 0$. It is a standard fact that $M$ is torsion over $\Zp\ps{G}$ if and only if $\Hom_{\Zp\ps{G}}(M,\Zp\ps{G})=0$ (for instance, see \cite[Lemma 4.2]{LimFine}). In the event that the torsion $\Zp\ps{G}$-module $M$ satisfies $\Ext^1_{\Zp\ps{G}}(M,\Zp\ps{G})=0$, we shall say that $M$ is a pseudo-null $\Zp\ps{G}$-module.

We now extend the notion of torsion modules and pseudo-null modules to the case when $G$ is a compact $p$-adic Lie group which is not necessarily pro-$p$. A well-known theorem of Lazard asserts that the $p$-adic Lie group $G$ contains an open normal subgroup $G_0$ which is pro-$p$ with no $p$-torsion (cf.\ \cite[Theorem 8.32]{DDMS}). By \cite[Proposition 5.4.17]{NSW}, we have
\[\Ext^i_{\Zp\ps{G}}(M,\Zp\ps{G}) \cong \Ext^i_{\Zp\ps{G_0}}(M,\Zp\ps{G_0})\]
for every finitely generated $\Zp\ps{G}$-module $M$. In view of this, we shall say that $M$ is a torsion $\Zp\ps{G}$-module (resp., psuedo-null $\Zp\ps{G}$-module) if $\Hom_{\Zp\ps{G}}(M,\Zp\ps{G})=0$ (resp., $\Ext^i_{\Zp\ps{G}}(M,\Zp\ps{G})=0$ for $i=0,1$). Equivalently, this is saying that $M$ is a torsion $\Zp\ps{G}$-module (resp., pseudo-null $\Zp\ps{G}$-module), whenever $M$ is a torsion $\Zp\ps{G_0}$-module (resp., pseudo-null $\Zp\ps{G_0}$-module) as in the preceding paragraph. (Also, compare with \cite[Discussion after Definition 2.6]{V02}).

\subsection{Characteristic elements in Iwasawa theory} \label{char element subsec}

Following \cite{Bu15, BV, CF+, FK}, we now describe how to attach characteristic elements to a certain class of $\Zp\ps{G}$-modules.
From now on, $G$ will always denote a compact $p$-adic Lie group which contains a closed normal subgroup $H$ such that $\Ga:= G/H \cong \Zp$. As before, we also assume that the group $G$ has no $p$-torsion.

\bd
Denote by $\M_H(G)$ the category of finitely generated $\Zp\ps{G}$-module $M$ with the property that $M/M[p^\infty]$ is finitely generated over $\Zp\ps{H}$.

This important class of modules was first introduced and studied intensively in \cite{CF+}. It has further been conjectured
that the dual Selmer group of an $p$-ordinary elliptic curve lies in the category $\M_H(G)$ (see \cite{BV, CF+, FK}; also see Conjecture \ref{MHG conj} below, and the discussion that follows concerning about when the said conjecture is known to hold).
\ed

To facilitate subsequent discussion, we recall an equivalent description of this class of modules as given in \cite{CF+}. Set
\[\Si : = \Si_{G,H} := \big\{ x\in \Zp\ps{G}~ \big|~  \Zp\ps{G}/\Zp\ps{G}x \mbox{ is finitely generated over } \Zp\ps{H} \big\}.\]
Thanks to \cite[Theorem 2.4]{CF+}, we now know that $\Si$ is a left and right Ore set of $\Zp\ps{G}$. Setting $\Si^* = \cup_{n\geq 0}p^n\Si$, it therefore makes sense to speak of the localization of $\Zp\ps{G}$ with respect to $\Si^*$, and the resulting ring is in turn denoted by  $\Zp\ps{G}_{\Si^*}$. Furthermore, it follows from \cite[Proposition 2.3]{CF+} that a finitely generated $\Zp\ps{G}$-module $M$ is annihilated by $\Si^*$ if and only if $M/M[p^\infty]$ is finitely generated over $\Zp\ps{H}$. It follows from this that there is an identification $K_0(\Zp\ps{G}, \Zp\ps{G}_{\Si^*}) \cong K_0(\M_H(G))$ (see \cite[Section 1.1]{BV}).

We now come to the process of attaching characteristic elements to modules in $\M_H(G)$. For this, we recall that the localization sequence in $K$-theory yields
the following exact sequence
\[ K_1(\Zp\ps{G}) \lra K_1(\Zp\ps{G}_{\Si^*})\stackrel{\partial_G}{\lra} K_0(\Zp\ps{G}, \Zp\ps{G}_{\Si^*}) \lra K_0(\Zp\ps{G}) \lra K_0(\Zp\ps{G}_{\Si^*})\lra 0\]
(for instance, see \cite{BerK}). The following result is fundamental for our subsequent discussion.

\bl \label{surjective connecting}
The above connecting homomorphism $\partial_G$ is surjective.
\el

\bpf
This was established by Coates \textit{et al} in \cite[Proposition 3.4]{CF+} (also see \cite[Corollary 3.8]{Wi}).
\epf

By virtue of Lemma \ref{surjective connecting}, one can make the following definition.

\bd[Coates \textit{et al} \cite{CF+, V05}] \label{CoatesFKSV} Let $M$ be a module in $\M_H(G)$. Then a characteristic element for $M$ is defined to be any element $\xi_M\in K_1(\Zp\ps{G}_{\Si^*})$
 such that $\partial_G(\xi_M) = -[M]$.
\ed

\br
In the original definition of \cite{CF+, V05}, they have chosen the characteristic element such that $\partial_G(\xi_M) = [M]$. Here we have adopted the later convention in \cite{BV, FK} by having the characteristic element being sent to $-[M]$ (see especially \cite[Appendix C]{BV} for the rationale of this choice).
\er

We end the subsection describing how characteristic elements behave under restriction in $K$-theory.
Let $U$ be an open normal subgroup of $G$. Set $H_U:= H\cap U$ and $\Ga_U:= HU/H\cong U/H_U$. Write $\ga_U = \ga^{|\Ga:\Ga_U|}$ which is now a topological generator of the group $\Ga_U$. Plainly $\Zp\ps{G}$ is finite free over $\Zp\ps{U}$. Furthermore, if we write $\Si^*_U$ for $\Si^*_{U,H_U}$, then $\Zp\ps{G}_{\Si^*}$ is finite free over $\Zp\ps{U}_{\Si^*_U}$ (cf.\ \cite[Proposition 4.5(i)]{SchV10}). This in turn induces a natural restriction homomorphism
\[\res : K_1(\Zp\ps{G}_{\Si^*}) \lra K_1(\Zp\ps{U}_{\Si^*_U})\]
on the $K_1$-groups. Now, if $M$ lies in $\M_H(G)$, then it is straightforward to verify that $M$ lies in $\M_{H_U}(U)$ too. Therefore, it makes sense to speak of a characteristic element of $M$ which now lies in $ K_1(\Zp\ps{U}_{\Si^*_U})$.

\bl \label{res BVseries}
 Retain notation as above. If $\xi_M$ is a characteristic element of $M$, then $\res(\xi_M)$ is a characteristic element of $M$ when viewed in  $\M_{H_U}(U)$.
\el

\bpf
By the functoriality of $K$-groups, we have the following commutative diagram
\[  \xymatrixcolsep{0.6in}
\entrymodifiers={!! <0pt, .8ex>+} \SelectTips{eu}{}\xymatrix{
    K_1(\Zp\ps{G}_{\Si^*})  \ar[r]^{\partial_G} \ar[d]^{\res} & K_0(\M_H(G))  \ar[d]^{\res}  \\
    K_1(\Zp\ps{U}_{\Si^*_U}) \ar[r]^{\partial_U} &  K_0(\M_{H_U}(U)) } \]
Therefore, if $\xi_M$ is a characteristic element of $M\in\M_H(G)$, then
\[ \partial_U(\res(\xi_M)) = \res(\partial_G(\xi_M)) = res(-[M]) = -[M]\]
which is precisely saying that $\res(\xi_M)$ is a characteristic element of $M$ in the terminology of Definition \ref{CoatesFKSV}.
\epf

\subsection{Evaluation at Artin representations}

To continue, we shall fix an algebraic closure $\overline{\Q}_p$ of $\Qp$. Let $\rho :G \lra \mathrm{GL}_n(\Op)$ be an Artin representation, by which we mean that $\rho$ is a continuous group homomorphism with an open kernel, and where $\Op$ is the ring of integers of a finite extension of $\Qp$ contained in our fixed choice of $\overline{\Q}_p$.
For each $g \in G$, write $\bar{g}$ for its image in $\Ga= G/H$. We then define a group homomorphism
\[ G \lra \mathrm{GL}_n(\Op) \ot \Zp\ps{\Ga},\quad g \mapsto \rho(g)\ot\bar{g}. \]
By \cite[Lemma 3.3]{CF+}, this extends to a ring homomorphism
\[ \Zp\ps{G}_{\Si^*} \lra M_n(\Op) \ot Q(\Ga) \cong M_n(Q_{\Op}(\Ga)),  \]
where $Q(\Ga)$ (resp., $Q_{\Op}(\Ga)$) denotes the field of fraction of $\Zp\ps{\Ga}$ (resp., $\Op\ps{\Ga}$). It then follows from the functoriality of $K_1$-groups that we have a group homomorphism $ K_1(\Zp\ps{G}_{\Si^*}) \lra  K_1(M_n(Q_{\Op}(\ps{\Ga})))$ which in turn fits into the following composition
\[ \Phi_\rho: K_1(\Zp\ps{G}_{\Si^*}) \lra  K_1(M_n(Q_{\Op}(\ps{\Ga}))) \stackrel{\mathrm{Morita}}\cong K_1(Q_{\Op}(\ps{\Ga})) \cong Q_{\Op}(\ps{\Ga})^{\times} \cong Q_{\Op}(T)^{\times},\]
where the final isomorphism is given by $\ga-1 \mapsto T$ for a fixed topological generator $\ga$ of $\Ga$, and $Q_{\Op}(T)$ is the field of fraction of the power series ring $\Op\ps{T}$ in one variable. We now make a preliminary definition.

\bd[Burns \cite{Bu15}] \label{def ord vanishing}
Let $\xi\in K_1(\Zp\ps{G}_{\Si^*})$ and $\rho$ an Artin representation of $G$. Fix a topological generator $\ga$ of $\Ga$. Then one has $\Phi_{\rho}(\xi) = T^{r_{\rho,\ga}(\xi)}g_{\ga}(T)$ for some integer $r_{\rho,\ga}(\xi)$ and $g_{\ga}(T)\in Q_{\Op}(T)^{\times}$ such that $g_{\ga}(0)\neq 0$. 
In the event that the integer $r_{\rho,\ga}(\xi)$ is non-negative, we write $\xi^*(\rho,\ga)$ for the value $g_{\ga}(0)$. 
\ed

\br
In most arithmetic situations, especially for the case of Selmer groups, the non-negativity of the integers $r_{\rho,\ga}(\xi)$ is a consequence of a conjecture of Coates \textit{et al} (see \cite[Conjecture 4.8 and discussion after]{CF+}). Although this said conjecture is only known for a $\Zp^d$-extension (see \cite[p. 197]{CF+}) and a False-Tate extension \cite{K05}, the non-negativity of the integer $r_{\rho,\ga}(\xi)$ can be verified by a direct calculation of the so-called Akashi series (see Subsection \ref{Ak subsec}, and Sections \ref{main results} and \ref{twisted Section}). This thankfully suffices for our discussion in the paper.
\er

At first viewing, the quantities in Definition \ref{def ord vanishing} seem to rely on the choice of the generator of $\Ga$. We shall see immediately that this is not the case.

\bl \label{independence of generator and units}
Retain notations as in Definition \ref{def ord vanishing}. Then the following statements hold.
\begin{enumerate}
\item[$(i)$] The integer $r_{\rho,\ga}(\xi)$ is independent of the choice of generator $\ga$ of $\Ga$. In the event that $r_{\rho,\ga}(\xi)$ is non-negative, the value $\xi^*(\rho,\ga)$ mod $\Op^\times$ is also independent of the choice of $\ga$.
\item[$(ii)$]  Suppose that $M$ is a module belonging to $\M_H(G)$. Then the quantities $r_{\rho,\ga}\big(\xi_M)$ and $\xi_M^*(\rho,\ga)$ mod $\Op^\times$ are independent of the choice of the characteristic element $\xi_M$.
\end{enumerate}

\el

\bpf
(i) Let $\ga'$ be another generator of $\Ga$. Then we have $\ga' = \ga^u$ for some $u\in \Zp^{\times}$. Identify $Q_{\Op}(\Ga)^{\times}\cong Q_{\Op}(T')^{\times}$ via $\ga'\mapsto T'+1$. With respect to $\ga'$, we have $\Phi_{\rho}(\xi) = (T')^r g(T')$ for some integer $r$ and $g(T')\in Q_{\Op}(T')^{\times}$ such that $g(0)\neq 0$. Since $\ga' = \ga^u$, we have $T' = (T+1)^u-1$ and performing this substitution, we obtain
\[\Phi_{\rho}(\xi) = \big((T+1)^u-1\big)^r g((T+1)^u-1) \]
in $Q_{\Op}(T)^{\times}$. Direct calculations show that $g\big((T+1)^u-1\big)\big|_{T=0} = g(0) \neq 0$ and
\[\frac{\big((T+1)^u-1\big)^r}{T^r}\Big|_{T=0} = u^r\in \Zp^{\times}\subseteq \Op^{\times}.\]
In other words, with respect to $\ga$, we have $\Phi_{\rho}(\xi) = T^r h(T)$ for some $h(T)\in Q_{\Op}(T)$ such that $h(0) = g(0)u^r$  mod $\Op^\times$. Assertion (i) now follows.

(ii) Two characteristic elements of $M$ must differ by an element $z$ which lies in the image of the map \[K_1(\Zp\ps{G}) \lra K_1(\Zp\ps{G}_{\Si^*}).\]
 But by \cite[Lemma 4.9]{CF+}, we have $\Phi_{\rho}(z) \in \Op^\times$. This proves assertion (ii).
\epf

In view of the preceding lemma, we can now make the following definition.

\bd \label{BurnsDef}
Let $\xi\in K_1(\Zp\ps{G}_{\Si^*})$ and $\rho$ an Artin representation of $G$. We write
\[r_{\rho}(\xi)= \ord_{T=0} \big(\Phi_{\rho}(\xi)\big)\]
for the order of vanishing as defined in Definition \ref{def ord vanishing}.

In the event that this order of vanishing is non-negative, we write
\[ \frac{1}{T^{r_{\rho}(\xi)}}\Phi_{\rho}(\xi)\Big|_{T=0}\]
for $\xi^*(\rho)$ mod $\Op^\times$.
\ed

We end the subsection with a result on the behavior of the quantities defined in Definition \ref{BurnsDef} under induction of characters. This will play an important role in the eventual proof of our main result.
Let $U$ be a given open normal subgroup of $G$. We set $H_U:= H\cap U$ and $\Ga_U:= U/H_U$. Write $\ga_U = \ga^{|\Ga:\Ga_U|}$ which is a topological generator of $\Ga_U$. Recall that there is a natural restriction homomorphism
\[\res : K_1(\Zp\ps{G}_{\Si^*}) \lra K_1(\Zp\ps{U}_{\Si^*_U}). \]
For each $\eta \in K_1(\Zp\ps{U}_{\Si^*_U})$ and an Artin representation $\psi$ of $U$, we write $\Phi_{\psi,U}(\eta)$ for the evaluation of $\eta$ at $\psi$ which lies in $Q_{\Op}(T_U)^{\times}$, where $Q_{\Op}(T_U)^{\times}\cong Q_{\Op}(\Ga_U)^{\times}$ via $T_U\mapsto \ga_U-1$. Under this choice of identification, we may view $Q_{\Op}(T_U)^{\times}\subseteq Q_{\Op}(T)^{\times}$, where $T_U$ is sent to $(1+T)^{d_U}-1$ with $d_U = |\Ga:\Ga_U|$.

\bp \label{evaluation ind rep}
Let $U$ be an open normal subgroup of $G$ and $\psi$ an Artin representation of $U$. Set $\rho = \mathrm{Ind}^G_U\psi$. Then the following statements are valid.
\begin{enumerate}
  \item[$(i)$] $\Phi_{\psi, U}(\res ~\xi) = \Phi_{\rho, G}(\xi)$ in $Q_{\Op}(T_U)^{\times}\subseteq Q_{\Op}(T)^{\times}$.
  \item[$(ii)$] $\ord_{T_U=0}\big(\Phi_{\psi, U}(\res ~\xi)\big) = \ord_{T=0}\big(\Phi_{\rho, G}(\xi)\big)$.
  \item[$(iii)$] $(d_U)^r\res (\xi)^*(\psi) = \xi^*(\rho)~ \mbox{mod} ~\Op^\times$, where $r$ is the common value in $(ii)$.
\end{enumerate}
\ep

\bpf
 This is essentially \cite[Lemma 3.6]{Bu15} and we sketch the idea behind it. We should mention that the first identity is quite subtle. In fact, it can be presented in the form of the following commutative diagram.
 \[  \xymatrixcolsep{0.6in}
\entrymodifiers={!! <0pt, .8ex>+} \SelectTips{eu}{}\xymatrix{
    K_1(\Zp\ps{G}_{\Si^*})  \ar[r]^{\Phi_{\rho, G}} \ar[d]^{\res} & Q_{\Op}(T)^{\times}    \\
    K_1(\Zp\ps{U}_{\Si^*_U}) \ar[r]^{\Phi_{\psi, U}} & Q_{\Op}(T_U)^{\times} \ar @{^{(}->}[u] } \]
(Take note of the direction of the rightmost vertical arrow!) One first reduces the problem to the case of a dimension one group $G$. When in the one-dimensional case, the evaluation map $\Phi_{\rho, G}$ can be related to the reduced norm in the sense of Ritter and Weiss \cite{RW} (see \cite[Lemma 3.1]{Bu15}). Under this identification, the commutativity of the above diagram then follows from a calculation of Ritter-Weiss \cite[Lemma 9]{RW}.

Now, by virtue of (i), one has
 \[ r: = \ord_{T_U=0}\big(\Phi_{\psi, U}(\res ~\xi)\big) = \ord_{T_U=0}\big(\Phi_{\rho, G}(\xi)\big) .\]
 In other words, we have $\Phi_{\rho, G}(\xi) = T_U^r g(T_U)$ with $g(T_U)|_{T_U=0}\neq 0$. Substituting $T_U$ for $(1+T)^{d_U}-1$, we obtain
 \[\Phi_{\rho, G}(\xi) = \big((1+T)^{d_U}-1\big)^r g\big((1+T)^{d_U}-1\big). \]
 A direct calculation shows that $g\big((T+1)^{d_U}-1\big)\big|_{T=0} = g(T_U)|_{T_U=0} \neq 0$ and
\[\frac{\big((T+1)^{d_U}-1\big)^r}{T^r}\Big|_{T=0} = (d_U)^r\neq 0.\]
This thus yields (ii) and (iii).
\epf

\subsection{Akashi series} \label{Ak subsec}

We continue to let $G$ denote a compact $p$-adic Lie group without $p$-torsion, which contains a closed normal subgroup $H$ such that $\Ga:= G/H \cong \Zp$. We now introduce an algebraic invariant that was first defined by Coates-Schneider-Sujatha in \cite{CSS} which has played an important role in the study of Selmer groups over non-commutative $p$-adic Lie extensions (see \cite{CF+, Ze11}).

\bd[\cite{CSS}]
Let $M$ be a finitely generated $\Zp\ps{G}$-module. Suppose that $M$ has the property that $H_i(H,M)$ is a torsion $\Zp\ps{\Ga}$-module for every $i$. We define the Akashi series of $M$ to be
\[ \Ak_H(M) = \prod_i \ch_{\Zp\ps{\Ga}}\big(H_i(H,M)\big)^{(-1)^i},\]
where $\ch_{\Zp\ps{\Ga}}(N)$ is the characteristic element of the $\Zp\ps{\Ga}$-module $N$. Note that the Akashi series is only well-defined up to a unit in $\Zp\ps{\Ga}$, although this suffices for most purposes. In the event that the Akashi series of $M$ is a
unit in $\Zp\ps{\Ga}$, we shall sometimes write $\Ak_H(M) = 1$. The following lemma gives a useful criterion for determining the unicity of the Akashi series.
\ed

\bl \label{Akashi Zp}
Suppose that $G$ is a pro-$p$ group containing a closed normal subgroup $N$ contained in $H$ such that $G/N= U\times H/N$, where $H/N$ is of dimension $\geq 1$ and $U\cong G/H$. For every $\Zp\ps{G}$-module $M$ that is finitely generated
over $\Zp$, we have $\Ak_H(M) = 1$.
\el

\bpf
See \cite[Proposition 2.3]{Ze11}.
\epf

For our purposes, we require the following basic result which relates the evaluation of $\xi_M$ for a module $M$ in $\M_H(G)$ to its twisted Akashi series. For a given Artin representation $\rho$, we let $W_\rho$ denote a free $\Op$-module of rank $d$ realizing $\rho$. Set $\mathrm{tw}_{\rho}(M) = M\ot_{\Zp}W_\rho$, which is endowed with a diagonal action of $G$, i.e., for $g\in G$, we have $g(m\ot w) = gm \ot \rho(g)w$. Finally, we write $\hat{\rho}$ for the contragradient of $\rho$, i.e., $\hat{\rho}(g) = \rho(g^{-1})^t$ for $g\in G$, where $t$ denotes the transpose matrix. With these notations in hand, we can now state the next lemma.

\bl \label{char and Akashi}
Suppose that the $p$-adic Lie group $G$ has no $p$-torsion.
Let $M$ be a module in $\M_H(G)$ and $\xi_M\in K_1(\Zp\ps{G}_{\Si^*})$ a characteristic element of $M$. Then we have
 $\Phi_\rho(\xi_M) = \Ak_H(\mathrm{tw}_{\hat{\rho}}(M))$ mod $\Op\ps{\Ga}^\times$.
\el

\bpf
By \cite[Lemmas 3.1 and 3.2]{CF+}, $H_i(H,\tw_{\hat{\rho}}(M))$ is a torsion $\Zp\ps{\Ga}$-module for every $i$. Therefore, the quantity $\Ak_H(\mathrm{tw}_{\hat{\rho}}(M))$ is well-defined. The required formula now follows from \cite[Lemma 3.7]{CF+}.
\epf

\section{Elliptic curves over local fields} \label{Elliptic local field}

In this section, we record certain results on elliptic curves over a $p$-adic local field. As a start, we introduce a notation which will be adhered for the remainder of the paper without further mention. If $K$ is a field and $W$ is a $\Gal(\bar{K}/K)$-module, write $W(K) = W^{\Gal(\bar{K}/K)}$. In the event that $K$ is a local field of characteristic zero, we let $K^\cyc$  denote the cyclotomic $\Zp$-extension of $K$, and write $\Ga= \Gal(K^\cyc/K)\cong\Zp$.
We shall always identify $\Zp\ps{\Ga}$ with $\Zp\ps{T}$ under a fixed choice of topological generator of $\Ga$.

We consider the good ordinary reduction case first.

\subsection{Good ordinary reduction} \label{good ord subsec}

Let $E$ be an elliptic curve defined over a finite extension $K$ of $\Qp$ which has good ordinary reduction. Then from
\cite{CG} (or \cite{G99}), we have the following short exact sequence of $\Gal(\bar{K}/K)$-modules
\[ 0 \lra \widehat{E}_{p^\infty} \lra \Ep \lra  \widetilde{E}_{p^\infty}\lra 0,\]
where $\widehat{E}$ (resp., $\widetilde{E}$) is the formal group (resp., reduced curve) of $E$. Furthermore, the group $\widetilde{E}_{p^\infty}$ is $p$-divisible of $\Zp$-corank one with a trivial action of the inertia subgroup of $\Gal(\bar{K}/K)$. Furthermore, by \cite[Propositions 4.7 and 4.8]{CG}, if $K_{\infty}$ is a $p$-adic Lie extension of $K$ which is infinitely ramified, we have an isomorphism
\[ H^1(K_\infty, E)_{p^\infty} \cong H^1(K_\infty,\widetilde{E}_{p^\infty}). \]
This identification will be frequently utilized in the subsequent discussion of the paper. In the event that $K_\infty$ contains the cyclotomic $\Zp$-extension $K^\cyc$, the above identification can be applied to both $K^\cyc$ and $K_\infty$. In view of these identifications, the subsequent discussion in Sections 4-6 necessitates us to understand the $\Gal(K_\infty/K^\cyc)$-cohomology of $\widetilde{E}_{p^\infty}$. For this, we have the following.

\bl \label{akashi ord p}
 Suppose that $K_\infty$ is a pro-$p$ extension of $K$ which satisfies the following properties.
 \begin{enumerate}
 \item[$(a)$] The field $K_\infty$ contains the cyclotomic $\Zp$-extension $K^\cyc$.
 \item[$(b)$] The group $\Gal(K_\infty/K)$ has no $p$-torsion.
 \end{enumerate}
 Writing $\mH=\Gal(K_\infty/K^\cyc)$ and $\Ga=\Gal(K^\cyc/K)$, we have the following equalities.
 \begin{enumerate}
 \item[$(i)$] $\Ak_\mH(\widetilde{E}_{p^\infty}(K_\infty)^{\vee})=1$.
 \item[$(ii)$] $\displaystyle\prod_{j\geq 1}\ch_{\Zp\ps{\Ga}}\left( H^j(\mH, \widetilde{E}_{p^\infty}(K_\infty))^{\vee}\right)^{(-1)^{j+1}} =1$.
 \end{enumerate}
\el

\bpf
 Recall that we have
    \[\Ak_\mH(\widetilde{E}_{p^\infty}(K_\infty)^{\vee}) = \ch_{\Zp\ps{\Ga}}\big(H^0(\mH, \widetilde{E}_{p^\infty}(K_\infty))^\vee\big)\times \displaystyle\prod_{j\geq 1}\ch_{\Zp\ps{\Ga}}\left( H^j(\mH, \widetilde{E}_{p^\infty}(K_\infty))^{\vee}\right)^{(-1)^{j}}.\]
 Since $K^\cyc/K$ is a ramified $\Zp$-extension, $\widetilde{E}_{p^\infty}(K^\cyc)= H^0(\mH, \widetilde{E}_{p^\infty}(K_\infty))$ must be finite and so has trivial $\Zp\ps{\Ga}$-characteristics. Therefore, the equality in $(ii)$ will follow once we establish $(i)$ which we will do. Now, since $\widetilde{E}_{p^\infty}$ is divisible of $\Zp$-corank one, it is realizable over $K_\infty$ if and only if $\widetilde{E}_{p^\infty}(K_\infty)$ is infinite. Therefore, if $\widetilde{E}_{p^\infty}$ is not realizable over $K_\infty$, then $\widetilde{E}_{p^\infty}(K_\infty)$ is finite and so $H^j(\mH, \widetilde{E}_{p^\infty}(K_\infty)$ is finite for every $j$. In particular, this implies that $\Ak_\mH(\widetilde{E}_{p^\infty}(K_\infty)^{\vee})=1$.

 Now suppose that $\widetilde{E}_{p^\infty}(K_\infty) = \widetilde{E}_{p^\infty}$. Then we have that $K^\cyc(\widetilde{E}_{p^\infty})$ is contained in $K_\infty$. Since $\Gal(K_\infty/K)$ is pro-$p$ with no $p$-torsion, we must have $\Gal(K^\cyc(\widetilde{E}_{p^\infty})/K) \cong \Zp^2$. The unicity of the Akashi series now follows from an application of Lemma \ref{Akashi Zp}.
\epf

We also consider a twisted variant of the preceding lemma. Recall that if $\rho$ is an Artin representation of $\Gal(K_\infty/K)$ with coefficient in $\Op$, we write $\tw_{\rho}(-)= -\ot_{\Zp}W_{\rho}$, where $W_\rho$ denotes the free $\Op$-module realizing the representation $\rho$.

\bl \label{akashi ord p twist}
 Suppose that $K_\infty$ is a $p$-adic extension of $K$ which satisfies all of the following properties.
 \begin{enumerate}
 \item[$(a)$] The field $K_\infty$ contains the cyclotomic $\Zp$-extension $K^\cyc$.
 \item[$(b)$] The group $\Gal(K_\infty/K)$ has no $p$-torsion $($but not necessarily pro-$p$$)$.
 \item[$(c)$] $H^j(\mathcal{H}', \widetilde{E}_{p^\infty}(K_\infty))$ is finite for every $j\geq 1$ and open subgroup $\mH'$ of $\mH=\Gal(K_\infty/K^\cyc)$.
 \end{enumerate}
 Then $H^j\big(\mH, \tw_{\rho}(\widetilde{E}_{p^\infty})(K_\infty)\big)$ is finite for every $j\geq 1$.
  \el

\bpf
  Let $\mathcal{H}'$ be an open subgroup of $\mH$ such that $\mH'\subseteq \ker \rho $ (for instance, take $\mathcal{H}'=\mathcal{H}\cap\ker\rho$). In view of the following spectral sequence
  \[ H^r\big(\mH/\mH', H^s(\mH', \tw_{\rho}(\widetilde{E}_{p^\infty})(K_\infty)) \big) \Longrightarrow H^{r+s}\big(\mH, \tw_{\rho}(\widetilde{E}_{p^\infty})(K_\infty)\big), \]
it suffices to show that the term $E_{r,s}:=H^r\big(\mH/\mH', H^s(\mH', \tw_{\rho}(\widetilde{E}_{p^\infty})(K_\infty)) \big)$ is finite for every $r, s$.
  By our choice of $\mH'$, whenever $s\geq 1$, we have isomorphisms
 \[ H^s\big(\mH', \tw_{\rho}(\widetilde{E}_{p^\infty})(K_\infty) \big) \cong H^s\big(\mH',  \widetilde{E}_{p^\infty}(K_\infty) \big)\ot_{\Zp}W_{\rho}  \]
 of abelian groups. Taking hypothesis (c) into account, this yields the finiteness of $E_{r,s}$ for $s\geq 1$.

  Now suppose that $s=0$. Let $L'_\infty$ denote the fixed field of $K_\infty$ by $\mH'$. Then one has
 \[H^r\big(\mH/\mH', H^0(\mH', \tw_{\rho}(\widetilde{E}_{p^\infty})(K_\infty)) \big) = H^r\big(\mH/\mH', \tw_{\rho}(\widetilde{E}_{p^\infty})(L'_\infty)\big).  \]
 On the other hand, it follows from our choice of $\mH'$ that
 \[\tw_{\rho}(\widetilde{E}_{p^\infty})(L'_\infty) = H^0\big(L'_\infty, \tw_{\rho}(\widetilde{E}_{p^\infty})\big) = H^0\big(L'_\infty, \widetilde{E}_{p^\infty}\big)\ot_{\Zp}W_{\rho}. \]
 As $\mH'$ is a subgroup of $\mH$ of finite index, the field $L'_\infty$ is a finite extension of $K^\cyc$. In particular, $L'_\infty$ is a ramified $\Zp$-extension of some finite extension of $K$. Hence
 $H^0\big(L'_\infty, \widetilde{E}_{p^\infty}\big)$ has to be finite, and this completes the proof of the lemma.
\epf

\subsection{Split multiplicative case}

In this subsection, we let $E$ denote an elliptic curve defined over $\Qp$ which has split multiplicative reduction. Therefore, there is a short exact sequence of $\Gal(\bar{\Q}_p/\Qp)$-modules
\[ 0 \lra \mu_{p^{\infty}} \lra \Ep \lra  \Qp/\Zp\lra 0\]
(cf. \cite[pp. 69-70]{G99}).

Let $K$ be a finite extension of $\Qp$. By restriction of scalars, the above can also be viewed as a short exact sequence of $\Gal(\bar{\Q}_p/K)$-modules which in turn induces the following exact sequence
\[ H^1(K, \mu_{p^{\infty}}) \stackrel{\la}{\lra} H^1(K, \Ep) \lra H^1(K, \Qp/\Zp) \lra H^2(K,  \mu_{p^{\infty}}) \lra 0, \]
where the rightmost zero follows from the fact that $H^2(K,\Ep)=0$ (see \cite[Proof of Lemma 1.12]{CS}).
Since $H^1(K, \mu_{p^{\infty}})$ is $p$-divisible, so is its image under $\la$. In view of this, we may apply \cite[Proposition 4.5]{CG} to conclude that $\im \la = \im \kappa$, where
\[\kappa : E(K)\ot\Qp/\Zp\lra H^1(K, \Ep)\]
is the Kummer homomorphism. From which, we obtain the following short exact sequence
\[ 0\lra H^1(K, E)_{p^{\infty}} \lra H^1(K, \Qp/\Zp) \lra H^2(K,  \mu_{p^{\infty}}) \lra 0. \]

For every finite extension of $K$ contained in $K^\cyc$, we have a similar short exact sequence as above. Upon taking direct limit, we obtain
\[H^1(K^\cyc, E)_{p^{\infty}} \cong H^1(K^{\cyc}, \Qp/\Zp),\]
noting that $H^2(K^\cyc, \mu_{p^\infty})=0$ by \cite[Theorem 7.1.8(i)]{NSW}. All of these fit into the following commutative diagram
\[
 \entrymodifiers={!! <0pt, .8ex>+} \SelectTips{eu}{}\xymatrix{
    0 \ar[r]^{} & H^1(K, E)_{p^{\infty}} \ar[d]^{r} \ar[r] &
    H^1(K, \Qp/\Zp)
    \ar[d]^{h} \ar[r] & H^2(K, \mu_{p^{\infty}}) \ar[r] & 0 \\
   & \big(H^1(K^{\cyc}, E)_{p^{\infty}}\big)^{\Ga} \ar[r]^{\cong} &  H^1(K^\cyc, \Qp/\Zp)^{\Ga}  & &
     } \]
with exact rows, where $\Ga = \Gal(K^\cyc/K)$. We can now state the main result of this subsection.

\bl \label{split mult finite}
Keep notations as above. The map $r$ has a finite kernel.
\el

\bpf
Let $q_E$ be the Tate period of $E$. Since $E$ is defined over $\Qp$, the Tate period $q_E$ lies in $\Qp$. By the theorem of Barr\'e-Sirieix, Diaz, Gramain and Philibert \cite{BSDFP}, $q_E$ is transcendental and so $\log_p(q_E)\neq 0$. Thus, writing $N_{K/\Qp}$ for the norm map, we have
\[\log_p(N_{K/\Qp}(q_E)) = \log_p(q_E^{|K:\Qp|}) = |K:\Qp|\log_p(q_E) \neq 0.\]
By \cite[discussion in p. 78]{G99}, this in turn implies that the map $r$ has a finite kernel.
\epf

In this split multiplicative context, the results of Coates-Greenberg \cite[Propositions 4.7 and 4.8]{CG} yield an identification
\[ H^1(K_\infty, E)_{p^\infty} \cong H^1(K_\infty,\Qp/\Zp). \]
for an infinitely ramified $p$-adic Lie extension $K_{\infty}$ of $K$. Similarly to that in Subsection \ref{good ord subsec}, we are therefore led to studying the $\Gal(K_\infty/K^\cyc)$-cohomology of $\Qp/\Zp$. In this aspect, we shall prove two lemmas describing the order of vanishing of the characteristic elements of the $\Gal(K_\infty/K^\cyc)$-cohomology of $\Qp/\Zp$ in certain classes of $p$-adic Lie extensions of the local field $K$.  As before, we shall write $\Ga=\Gal(K^\cyc/K)$ and identify $\Zp\ps{\Ga}$ with $\Zp\ps{T}$ under a fixed choice of a topological generator of $\Ga$.

\bl \label{akashi split mult 1}
 Suppose that $K_\infty$ is a pro-$p$ extension of $K$ which satisfies the following properties.
 \begin{enumerate}
 \item[$(a)$] The field $K_\infty$ contains the cyclotomic $\Zp$-extension $K^\cyc$.
 \item[$(b)$] The group $\Gal(K_\infty/K)$ has no $p$-torsion.
 \item[$(c)$] There exists a subextension $\mL_\infty$ of $\mK_\infty$ containing $K^\cyc$ with the property that $\Gal(\mL_\infty/K) \cong \Gal(\mL_\infty/K^\cyc)\times \Gal(K^\cyc/K)$, where $\Gal(\mL_\infty/K^\cyc)$ is a pro-$p$ group of dimension $\geq 1$.
 \end{enumerate}
 Writing $\mH=\Gal(K_\infty/K^\cyc)$, we have the following equalities.
 \begin{enumerate}
 \item[$(i)$] $\Ak_\mH(\Zp)=1$.
 \item[$(ii)$] $\displaystyle\prod_{j\geq 1}\left(\ch_{\Zp\ps{\Ga}}\big( H^j(\mH, \Qp/\Zp)^{\vee}\big)\right)^{(-1)^{j+1}} =T$.
 \end{enumerate}
 \el

\bpf
 Since $\ch_{\Zp\ps{\Ga}}\left( H^0(\mH, \Zp)\right)= T$, it remains to verify the first equality. But this is immediate from assumption (c) and Lemma \ref{Akashi Zp}.
\epf

\bl \label{akashi zero order split mult 2}
 Suppose that $K$ contains a primitive $p$-th root of unity. For $d\geq 2$, we consider \[K_\infty = K\left(\mu_{p^\infty}, \sqrt[p^\infty]{\alpha_1} , \ldots, \sqrt[p^\infty]{\alpha_{d-1}}\right),\] where $\alpha_1, \ldots , \alpha_{d-1}\in K^{\times}$, whose image in $K^\times/(K^{\times})^p$ are linearly independent over $\Z/p\Z$.
 Writing $\mH=\Gal(K_\infty/K^\cyc)$, we have
 \[\ord_{T=0}\Big(\ch_{\Zp\ps{\Ga}}\big( H^j(\mH, \Qp/\Zp)^{\vee}\big)\Big) =0 \]
 for every $j\geq 1$. \el

\bpf
  We prove this by induction on $d$. Suppose that $d=2$. Then $H^j(\mH, \Qp/\Zp) = 0$ for $j\geq 2$ and
  \[ H^1(\mH, \Qp/\Zp)^{\vee} = \Hom(\Zp(1),\Qp/\Zp)^{\vee} = \Zp(1),\]
  which implies that
  \[\ord_{T=0}\Big(\ch_{\Zp\ps{\Ga}}\big( H^1(\mH, \Qp/\Zp)^{\vee}\big)\Big) =0. \]
  Now suppose that $d\geq 3$. Consider the subextension $L_\infty = K\left(\mu_{p^\infty}, \sqrt[p^\infty]{\alpha_1} , \ldots, \sqrt[p^\infty]{\alpha_{d-2}}\right)$, and write $N=\Gal(L_\infty/K)$. By our induction hypothesis, we have
  \[\ord_{T=0}\Big(\ch_{\Zp\ps{\Ga}}\big( H^j(N, \Qp/\Zp)^{\vee}\big)\Big) =0 \] for every $j\geq 1$. Let $Z=\Gal(K_\infty/L_\infty)\cong\Zp$. From the degeneration of the spectral sequence
  \[ H^i\big(Z, H^j(N, \Qp/\Zp)\big)\Longrightarrow H^{i+j}(\mH, \Qp/\Zp),\]
  we have short exact sequences
  \[ 0\lra H^1\big(Z, H^j(N, \Qp/\Zp)\big)\lra H^{j+1}(\mH, \Qp/\Zp) \lra H^0\big(Z, H^j(N, \Qp/\Zp)\big) \lra 0\]
  for $j\geq 0$. For $j\geq 1$, these sequences and our induction hypothesis immediately yield \[\ord_{T=0}\Big(\ch_{\Zp\ps{\Ga}}\big( H^{j+1}(N, \Qp/\Zp)^{\vee}\big)\Big) =0. \]
  When $j=0$, the short exact sequence reads as
\[ 0\lra H^1(Z, \Qp/\Zp)\lra H^{1}(\mH, \Qp/\Zp) \lra H^0\big(Z, H^1(N, \Qp/\Zp)\big) \lra 0.\]
  Again, by our induction hypothesis, the order of vanishing of the term on the right is zero. By identifying $Z$ with $\Gal\big(K^\cyc(\sqrt[p^\infty]{\alpha_{d-1}})/K^\cyc\big)$, we see that the order of vanishing of the leftmost term is also zero by the initial $d=2$ case.
  \epf

\subsection{Non-split multiplicative case}

We now consider an elliptic curve $E$ defined over a finite extension $K$ of $\Qp$ which has non-split multiplicative reduction. Then we have the following short exact sequence of $\Gal(\bar{K}_p/K)$-modules
\[ 0 \lra C \lra \Ep \lra  \Qp/\Zp\ot W_{\chi}\lra 0\]
(cf. \cite[pp. 69-70]{G99}), where $\chi$ is an unramified character of $\Gal(\bar{K}/K)$ which factors through a quadratic extension $K_\chi$ of $K$, and $W_\chi$ is the free $\Zp$-module (of rank 1) realizing $\chi$.
For later discussion, we record the following lemma.

\bl \label{Non-split local}
 Let $E$ be an elliptic curve defined over $K$ which has non-split multiplicative reduction. Let $K_\infty$ be a Galois extension of $K$ containing $K^\cyc$ such that $\Gal(K_\infty/K)$ is a $p$-adic Lie group of dimension $\geq 2$ with no $p$-torsion and that $K_\chi\cap K_\infty = K$. Let $\rho$ be an Artin representation of $\mG$. Write $\mH=\Gal(K_\infty/K^\cyc)$.
 Then $H^j\big(\mH, \tw_{\rho}(\Qp/\Zp\ot W_{\chi})(K_\infty)\big)$ is finite for every $j\geq 1$.
 \el

 \bpf
 Clearly, the lemma will follow once we show that $\tw_{\rho}(\Qp/\Zp\ot W_{\chi})(K_\infty)$ is finite. Note that
 \[\tw_{\rho}(\Qp/\Zp\ot W_{\chi})(K_\infty) = (\Qp/\Zp\ot W_{\chi})(K_\infty)\ot_{\Zp}W_\rho, \]
 and so we are reduced to establishing the finiteness of $(\Qp/\Zp\ot W_{\chi})(K_\infty)$. Since $\Qp/\Zp\ot W_{\chi}$ is divisible of $\Zp$-corank one, it suffices to show that it is not fixed under the action of $\Gal(\bar{K}/K_\infty)$. But this is an immediate consequence of the hypothesis that $K_\chi\cap K_\infty = K$.
 \epf

\br \label{nonsplit remark}
The hypothesis $K_\chi\cap K_\infty = K$ is equivalent to saying that $E$ has
non-split multiplicative reduction at every subextension of $K_\infty/K$ (see \cite[Chap.\ V, Lemma 5.2 and Theorem 5.3]{Sil}).
\er

\subsection{Elliptic curves over local fields of residue char $\neq p$}

In this subsection, $K$ is taken to be finite extension of $\Q_l$, where $l\neq p$.

\bl \label{neq l structure twist}
 Let $E$ be an elliptic curve defined over $K$, and let $K_\infty$ be a Galois extension of $K$ such that $\Gal(K_\infty/K)$ is a $p$-adic Lie group of dimension 2 with no $p$-torsion. Write $\mG=\Gal(K_\infty/K)$, $\mH=\Gal(K_\infty/K^\cyc)$ and $\Ga= \Gal(K^\cyc/K)$. Let $\rho$ be an Artin representation of $\mG$. Then
 \[ \ord_{T=0}\Big(\ch_{\Zp\ps{\Ga}}\big(H^1(\mH,W_\rho\ot_{\Zp}\Ep(K_\infty))^\vee\big)\Big) =0. \]
 \el

\bpf
It suffices to show that
\[\big(H^1(\mH,W_\rho\ot_{\Zp}\Ep(K_\infty))^\vee\big)_{\Ga}\] is finite, or equivalently,
 \[H^1\big(\mH,W_\rho\ot_{\Zp}\Ep(K_\infty)\big)^\Ga\]
 is finite. From the spectral sequence
\[H^i\big(\Ga, H^j(\mH, W_\rho\ot_{\Zp}\Ep(K_\infty))\big)\Longrightarrow H^{i+j}\big(\mG, W_\rho\ot_{\Zp}\Ep(K_\infty)\big),\]
we obtain a surjection
\[ H^1\big(\mG,W_\rho\ot_{\Zp}\Ep(K_\infty)\big) \tha H^1\big(\mH, W_\rho\ot_{\Zp}\Ep(K_\infty)\big)^\Ga. \]
It therefore remains to show that $H^1\big(\mG,W_\rho\ot_{\Zp}\Ep(K_\infty)\big)$ is finite.
Now, the low degree terms of the spectral sequence
\[H^i\big(\mG, H^j(K_\infty, W_\rho\ot_{\Zp}\Ep(K_\infty))\big)\Longrightarrow H^{i+j}\big(K, W_\rho\ot_{\Zp}\Ep(K_\infty)\big),\]
yields an exact sequence
\[ 0 \lra H^1(\mG,W_\rho\ot_{\Zp}\Ep(K_\infty)) \lra H^1(K,W_\rho\ot_{\Zp}\Ep) \lra H^1(K_\infty, W_\rho \ot_{\Zp}\Ep)^{\mG}.\]
A theorem of Iwasawa (cf. \cite[Theorem 7.5.3]{NSW}) tells us that $K_\infty$ has no non-trivial $p$-extension. Hence it follows that $H^1(K_\infty, W_\rho \ot_{\Zp}\Ep)=0$ and
\[H^1\big(\mG, W_\rho\ot_{\Zp}\Ep(K_\infty)\big) \cong H^1(K,W_\rho\ot_{\Zp}\Ep). \]
Therefore, the problem is reduced to establishing the finiteness of $H^1(K,W_\rho\ot_{\Zp}\Ep)$.
For this, we let $L$ be a finite Galois extension of $K$ contained in $K_\infty$ such that $\Gal(K_\infty/L)\subseteq \ker \rho$. Then we have the following exact sequence
\[ 0 \lra H^1(\Gal(L/K), W_\rho \ot_{\Zp}\Ep(L)) \lra H^1(K, W_\rho \ot_{\Op}\Ep) \lra H^1(L,W_\rho\ot_{\Zp}\Ep)^{\Gal(L/K)}. \]
Plainly, the leftmost term is finite. Therefore, to show that $H^1(K,W_\rho\ot_{\Zp}\Ep)$ is finite, it suffices to show that $H^1(L,W_\rho\ot_{\Zp}\Ep)$ is finite. But $H^1(L,W_\rho\ot_{\Zp}\Ep) = H^1(L, \Ep)\ot_{\Zp} W_\rho$ by  our choice of $L$, and so it remains to establish the finiteness of $H^1(L, \Ep)$. Since $\l\neq p$, this latter group is isomorphic to $H^1(L,E)_{p^\infty}$. By Tate-duality (cf.\ \cite[Chap.\ I, Corollary 3.4]{Mi}), this is isomorphic to $\Ep(L)^\vee$ which in turn is finite by Mattuck's theorem. This therefore completes the proof of the proposition.
\epf

\section{Selmer groups over $p$-adic Lie extensions} \label{Sel Sec}

We now come to arithmetic. Here, we fix the notation that we shall use throughout
the remainder of the paper. To start, we let $F$ be a number field and $E$ an elliptic curve defined over $F$. Let $S$ be a finite set of primes of $F$ which contains all the primes above $p$, the infinite primes and the primes of bad reduction of $E$. Denote by $F_S$ the maximal algebraic extension of $F$ which is unramified outside $S$. For each $v\in S$ and a finite extension $L$ of $F$, we set
\[ J_v(E/L) = \bigoplus_{w|v}H^1(L_w,E)_{p^\infty}.\]
If $\mathcal{L}$ is an infinite extension of $F$ contained in $F_S$, we define
\[ J_v(E/\mathcal{L}) = \ilim_L J_v(E/L),\]
where $L$ runs through all finite extensions  of $F$ contained in $\mathcal{L}$.

The classical ($p$-primary) Selmer group of $E$ over $\mathcal{L}$ is defined by
\[ \Sel(E/\mathcal{L})=\ker\Big(H^1(G_S(\mathcal{L}), \Ep)\lra \bigoplus_{v\in S}J_v(E/\mathcal{L})\Big),\]
where we write $G_S(\mathcal{L}) = \Gal(F_S/\mathcal{L})$. The Pontryagin dual of $\Sel(E/\mathcal{L})$ is then denoted by $X(E/\mathcal{L})$.

\medskip
\textbf{The following assumptions will be in full force for our elliptic curve $E$.}
\begin{itemize}
\item[$\mathbf{(S1)}$] The elliptic curve $E$ has either good ordinary reduction or multiplicative reduction at each prime of $F$ above $p$.

  \item[$\mathbf{(S2)}$] For each $v$ of $F$ above $p$ at which $E$ has split multiplicative reduction, we assume that there is an elliptic curve $E'_v$ over $\Qp$ with split multiplicative reduction such that $E/F_v$ is obtained from $E_v'$ via base changing from $\Qp$ to $F_v$.
\end{itemize}

In view of the above standing assumptions, for each prime $v$ of $F$ above $p$, we have a short exact sequence
\[ 0\lra C_v\lra \Ep\lra D_v\lra 0 \]
of discrete $\Gal(\bar{F}_v/F_v)$-modules, where
\[ D_v =\begin{cases} \widetilde{E_{v}}_{, p^\infty},  & \mbox{if $E$ has good ordinary reduction at $v$}, \\
\Qp/\Zp, & \mbox{if $E$ has split multiplicative reduction at $v$}, \\
\Qp/\Zp\ot_{\Zp} W_{\chi_v}, & \mbox{if $E$ has non-split multiplicative reduction at $v$.}\end{cases}\]
Here $W_{\chi_v}$ is the free $\Zp$-module of rank 1 realizing $\chi_v$, where $\chi_v$ is an unramified character of $\Gal(\bar{F_v}/F_v)$ which factors through a quadratic extension of $F_v$.

For our purposes, it is convenient to work with an equivalent description of the local terms $J_v(E/\mathcal{L})$, following an insight of Coates-Greenberg \cite{CG}.
Let $\mathcal{L}$ be an algebraic extension of $F$. For every non-archimedean prime $w$ of $\mathcal{L}$, write $\mathcal{L}_w$ for the union of the completions at $w$ of the finite extensions of $F$ contained in $\mathcal{L}$. If $w$ is a prime above $p$, we write $D_w = D_v$, where $v$ is a prime of $F$ below $w$. Finally, we shall always denote by $F^{\cyc}$ the cyclotomic $\Zp$-extension of $F$. With these in hand, we have the following lemma.

\bl \label{local coh description}
Let $\mathcal{L}$ be an algebraic extension of $F^\cyc$ which is unramified outside a set of finite primes of $F$. Then we have an isomorphism
\[ J_v(E/\mathcal{L})\cong \begin{cases} \ilim_{\mathcal{L}'}\bigoplus_{w|v}H^1(\mathcal{L}'_w, D_w),  & \mbox{if $v$ divides $p$}, \\
\ilim_{\mathcal{L}'}\bigoplus_{w|v}H^1(\mathcal{L}'_w, \Ep), & \mbox{if $v$ does not divide $p$},\end{cases} \]
where the direct limit is taken over all finite extensions $\mathcal{L}'$ of $F^\cyc$ contained in $\mathcal{L}$.
\el

\bpf
See \cite[Propositions 4.1, 4.7 and 4.8]{CG} or \cite[Lemma 4.1]{LimMHG}.
\epf

\subsection{Cyclotomic $\Zp$-extension}

We now specialize to the situation of a cyclotomic $\Zp$-extension.

\begin{conjecture}[Mazur, Schneider] \label{mazurconj}
Write $F^{\cyc}$ for the cyclotomic $\Zp$-extension of $F$. Then $X(E/F^{\cyc})$ is a torsion $\Zp\ps{\Ga}$-module, where $\Ga=\Gal(F^{\cyc}/F)$.
\end{conjecture}

The conjecture was first stated by Mazur in \cite{Maz} for elliptic curves that have good ordinary
reduction at all primes of $F$ above $p$. The form we stated here is a special case of that
in \cite{Sch85}. At present, the best result in support of the conjecture is
due to Kato \cite{K}, who has proven it when $F$ is abelian over $\Q$ and $E$ is an elliptic
curve defined over $\Q$ with ordinary reduction at $p$. A useful consequence of torsionness for us is the following.

\bp \label{cyclotomic torsion}
Let $E$ be an elliptic curve defined over $F$ which has either good ordinary reduction or split multiplicative reduction at each prime of $F$ above $p$. Then $X(E/F^\cyc)$ is torsion over $\Zp\ps{\Ga}$ if and only if $H^2(G_S(F^\cyc),\Ep)=0$ and there is a short exact sequence
\[ 0\lra  \Sel(E/F^\cyc)\lra H^1(G_S(F^\cyc),\Ep)\lra \bigoplus_{v\in S}J_v(E/F^\cyc)\lra 0.\]
\ep

\bpf
 See \cite[Proposition 3.3]{LimMHG}.
\epf

We continue to assume that $X(E/F^\cyc)$ is torsion over $\Zp\ps{\Ga}$. Then the structure theory of $\Zp\ps{\Ga}$-module tells us that there is a pseudo-isomorphism
\[X(E/F^\cyc) \sim \bigoplus_{i=1}^s\Zp\ps{\Ga}/p^{\al_i} \times \bigoplus_{j=1}^t\Zp\ps{\Ga}/f_j^{\be_j},\]
where each $f_j$ is irreducible in $\Zp\ps{\Ga}$ and is not an associate of $p$. With these notations, we can now state the following semi-simplicity conjecture of Greenberg \cite[Conjecture 1.12]{G99}.

\begin{conjecture}[Greenberg] \label{semisimple conj}
$\be_j =1$ for every $j$.
\end{conjecture}

The significance of this conjecture of Greenberg lies in the following observation (for instance, see \cite[Page 59]{G99}).

\bp \label{semisimple}
Let $E$ be an elliptic curve defined over $F$ which satisfies $\mathbf{(S1)}$ and $\mathbf{(S2)}$. Suppose that $X(E/F^\cyc)$ is torsion over $\Zp\ps{\Ga}$ and that Conjecture \ref{semisimple conj} is valid. Then we have
\[\ord_{T=0}\Big(\ch\big(X(E/F^\cyc)\big)\Big) = \corank_{\Zp}\big(\Sel(E/F)\big).\]
\ep

\bpf
 Consider the following commutative diagram
 \[   \entrymodifiers={!! <0pt, .8ex>+} \SelectTips{eu}{}\xymatrix{
    0 \ar[r] &\Sel(E/F)  \ar[d]^{\al} \ar[r] &  H^1\big(G_{S}(F),\Ep\big)
    \ar[d]^{h} \ar[r] & \displaystyle \bigoplus_{v\in S}J_v(E/F) \ar[d]^{\oplus g_v} \ar[r]& 0 \\
    0 \ar[r]^{} & \Sel(E/F^{\cyc})^{\Ga} \ar[r]^{} & \Big(H^1\big(G_{S}(F^{\cyc}),\Ep\big)\Big)^{\Ga} \ar[r] &
    \left(\displaystyle \bigoplus_{v\in S}J_v(E/F^{\cyc})\right)^{\Ga} &  } \]
    with exact rows. We shall show that $\al$ has finite kernel and cokernel. Since $\Ga$ has cohomological dimension 1, the map $h$ is surjective with kernel $H^1(\Ga, \Ep(F^\cyc))$. Since $\Ep(F^\cyc)$ is finite by a theorem of Ribet \cite{Ribet}, so is this kernel. Hence it remains to show that each map $g_v$ has finite kernel. For ordinary primes above $p$ and primes not
dividing $p$, this is established in the mist of proving the control theorem in the ordinary case (for instance, see \cite[Theorem 1.2 and Section 3]{G99}). In the event of a split multiplicative prime above $p$, the finiteness follows from Lemma \ref{split mult finite} (also see \cite[Proposition 3.7]{G99}) thanks to our hypothesis $\mathbf{(S2)}$. Finally, if $v$ is a non-split multiplicative prime of $E$ above $p$, we fix a prime of $F^\cyc$ above it and by abuse of notation denote it by $v$. Writing $\Ga_v$ for the decomposition group of $\Ga$ at $v$, we then have $\ker g_v = H^1(\Ga_v, \Qp/\Zp\ot W_{\chi_v}(F^\cyc_v))$. But this is finite for $\Qp/\Zp\ot W_{\chi_v}(F^\cyc_v)$ is finite as seen in the proof of Lemma \ref{Non-split local}.

In conclusion, the restriction map $\Sel(E/F)\stackrel{\al}{\lra} \Sel(E/F^{\cyc})^{\Ga}$ has finite kernel and cokernel (under the assumptions of our proposition). Taking Pontryagin dual, we see that the map $X(E/F^\cyc)_{\Ga}\lra X(E/F)$ has finite kernel and cokernel which in turn implies that
\[ \rank_{\Zp}\big(X(E/F^\cyc)_{\Ga}\big) = \corank_{\Zp}\big(\Sel(E/F)\big).\]
Finally, in view of the validity of Conjecture \ref{semisimple conj}, the left hand side of the equality is precisely \[\ord_{T=0}\Big(\ch\big(X(E/F^\cyc)\big)\Big).\]
The proof of the proposition is now complete.
\epf

To the best knowledge of the author, there seems very little evidence in literature on Conjecture \ref{semisimple conj}. The following is one simple criterion we know at present for proving in
some concrete examples. This will be useful for our discussion in Section \ref{examples section}.  We remark that our criterion is inspired by the discussion in \cite[Proposition 11.1]{Wu}. In the following lemma, if $M$ is a $\Zp\ps{\Ga}$-module, we write $\la(M)$ for its Iwasawa $\la$-invariant.

\bl \label{semisimple lemma}
Let $E$ be an elliptic curve defined over $F$ which satisfies $\mathbf{(S1)}$ and $\mathbf{(S2)}$. Suppose that $X(E/F^\cyc)$ is torsion over $\Zp\ps{\Ga}$ and that there exists a non-negative integer $n$ such that
\[\la\big(X(E/F^\cyc)\big) = \corank_{\Zp}\big(\Sel(E/F_n)\big),\]
where $F_n$ is the intermediate subextension of $F^\cyc/F$ with $|F_n:F|=p^n$.
Then Conjecture \ref{semisimple conj} is valid for $X(E/F^\cyc)$.
\el

\bpf
 Identify $\Zp\ps{\Ga}\cong\Zp\ps{T}$ and fix a pseudo-isomorphism
\[X(E/F^\cyc) \sim \bigoplus_{i=1}^s\Zp\ps{T}/p^{\al_i} \times \bigoplus_{k_0=1}^{a_0}\Zp\ps{T}/T^{\delta_{k_0}} \times \bigoplus_{k_1=1}^{a_1}\Zp\ps{T}/\Phi_{p}^{\delta_{k_1}}\times \cdots \times \bigoplus_{k_n=1}^{a_n}\Zp\ps{T}/\Phi_{p^n}^{\delta_{k_n}}\times \bigoplus_{j=1}^t\Zp\ps{T}/f_j^{\be_j},\]
where $\Phi_{p^i}(X)$ is the $p^i$-cyclotomic polynomial and $\Phi_{p^i} := \Phi_{p^i}(1+T)$, and where each $f_j$ is irreducible in $\Zp\ps{T}$ and is not an associate of $p$ or $\Phi_{p^i}$ $(1\leq i\leq n)$, and $\delta_{k_i}, \be_j\geq 1$.
Plainly, we have
\[  a_0 + (p-1)a_1+\cdots+ p^{n-1}(p-1)a_n= \rank_{\Zp}\big(X(E/F^\cyc)_{\Ga_n}\big).\]
On the other hand, as seen in the proof of Proposition \ref{semisimple}, the latter is equal to $\corank_{\Zp}\big(\Sel(E/F_n)\big).$
Hence we have
\begin{eqnarray*}
  a_0 + (p-1)a_1+\cdots +p^{n-1}(p-1)a_n\! &\leq& \!\sum_{k_0=1}^{a_0}\delta_{k_0}+(p-1)\sum_{k_1=1}^{a_1}\delta_{k_1} +\cdots + p^{n-1}(p-1)\sum_{k_n=1}^{a_n}\delta_{k_n} +\sum_{j=1}^t\be_j\deg(f_j) \\
  &=&\la\big(X(E/F^\cyc)\big) =\corank_{\Zp}\big(\Sel(E/F_n)\big) \\
  &=& a_0 + (p-1)a_1+\cdots+ p^{n-1}(p-1)a_n.
\end{eqnarray*}
This in turn forces $\delta_{k_i} =1$ for all $k_i$'s, and $t=0$. In particular, Conjecture \ref{semisimple conj} holds.
\epf

We give an example to illustrate Lemma \ref{semisimple lemma} (see Section \ref{examples section} for more examples, where the said lemma applies). Consider the elliptic curve $5692A1: y^2 = x^3 + x^2 -18x +25$ and take $p=3$. Write $\Q_n$ for the intermediate extension of $\Q$ contained $\Q^\cyc$ such that $|\Q_n:\Q|=3^n$. From \cite[Proposition 11.1]{Wu}, we have
 \[\rank_\Z(E(\Q))=2, \quad\rank_\Z(E(\Q_1))=6, \quad \rank_\Z(E(\Q_n))=12 ~\mbox{for}~n\geq 2, \]
and finiteness of $\sha(E/\Q^\cyc)[3^\infty]$. From this, we have
\[\la\big(X(E/F^\cyc)\big) = \corank_{\Zp}\big(\Sel(E/\Q_n)\big)\]
for $n\geq 2$. Lemma \ref{semisimple lemma} therefore applies. In particular, $X(E/\Q^\cyc)$ satisfies the semi-simplicity conjecture of Greenberg.
We now determine the structure of $X(E/\Q^\cyc)$ building on these data (compare with \cite[Proposition 11.1]{Wu}).
By the proof of Proposition \ref{semisimple}, we have
\[ \rank_{\Z_3}\big(X(E/\Q^\cyc)_{\Ga}\big) = \corank_{\Z_3}\big(\Sel(E/\Q)\big) =2.\]
Therefore, the $T$-primary part of $X(E/\Q^\cyc)$ must be pseudo-isomorphic to
\[ (\Z_3\ps{T}/T)^{\oplus 2}. \]
Similarly, since
\[ \rank_{\Z_3}\big(X(E/\Q^\cyc)_{\Ga_1}\big) = \corank_{\Z_3}\big(\Sel(E/\Q_1)\big) =6,\]
this forces $X(E/\Q^\cyc)^{\Ga_1}$
 to be pseudo-isomorphic to
\[ (\Z_3\ps{T}/T)^{\oplus 2}\times
(\Z_3\ps{T}/\Phi_3)^{\oplus 2}. \]
Finally, the equality
\[ \rank_{\Z_3}\big(X(E/\Q^\cyc)_{\Ga_2}\big) = 12\] tells us that the remaining factor in the characteristic polynomial of $X(E/\Q^\cyc)$ is $\Phi_9$. In conclusion, we have
 \[  X(E/\Q^\cyc) \sim (\Z_3\ps{T}/T)^{\oplus 2}\times
(\Z_3\ps{T}/\Phi_3)^{\oplus 2} \times \Z_3\ps{T}/\Phi_9. \]

\subsection{$p$-adic Lie extension} \label{p-adic Lie subsection}

We say that $F_\infty$ is a strongly admissible $p$-adic Lie extension of $F$ if $F_\infty$ is a Galois extension of $F$ which satisfies all of the following properties.

\begin{enumerate}
\item[(a)]  $\Gal(F_\infty/F)$ is a $p$-adic Lie group with no $p$-torsion.
\item[(b)] $F_\infty$ contains $F^\cyc$.
\item[(c)]  $F_\infty$ is unramified outside a finite set of primes.
\end{enumerate}

In the event that $\Gal(F_\infty/F)$ is pro-$p$, we shall call $F_\infty$ a strongly admissible pro-$p$ $p$-adic Lie extension of $F$.
Write $G=\Gal(F_\infty/F)$, $H=\Gal(F_\infty/F^\cyc)$ and $\Ga=\Gal(F^\cyc/F)$. For subsequent discussion, we enlarge our set $S$ of primes to also contain the ramified primes of $F_\infty/F$. We now state the following natural extension of Conjecture \ref{mazurconj} and its consequence.

\begin{conjecture} \label{mazurconj2}
Let $F_\infty$ be a strongly admissible $p$-adic Lie extension of $F$. Then $X(E/F_\infty)$ is torsion over $\Zp\ps{G}$.
\end{conjecture}

\bp \label{p-adic torsion}
Let $E$ be an elliptic curve defined over $F$ which has either good ordinary reduction or multiplicative reduction at each prime above $p$. Suppose that $F_\infty$ is a strongly admissible $p$-adic Lie extension of $F$ such that $X(E/F_\infty)$ is torsion over $\Zp\ps{G}$. Then the following assertions are valid.

\begin{enumerate}
\item[$(a)$] $H^2(G_S(F_\infty),\Ep)=0$.
\item[$(b)$] There is a short exact sequence
\[ 0\lra  \Sel(E/F_\infty)\lra H^1(G_S(F_\infty),E[p^\infty])\lra \bigoplus_{v\in S}J_v(E/F_\infty)\lra 0.\]
\end{enumerate}
\ep

\bpf
 See \cite[Proposition 3.3]{LimMHG}.
\epf

At our current knowledge, the torsionness of $X(E/F_\infty)$ is not enough for us to attach a characteristic element to it (see \cite{CSSAlg}). To circumvent this difficulty, Venjakob came up with an algebraic $K$-theoretical approach (see \cite{V05}; also see \cite{CF+}) to define a characteristic element. This however comes at the expense of the following $\M_H(G)$-conjecture \cite{CF+, CS12, DL, Lee}.

\begin{conjecture}\label{MHG conj}
The module $X(E/F_\infty)$ lies in the category $\M_H(G)$. In other words, $X_f(E/F_\infty):=X(E/F_\infty)/X(E/F_\infty)[p^\infty]$ is finitely generated over $\Zp\ps{H}$.
\end{conjecture}

This important conjecture is essential in the development of non-commutative Iwasawa theory, as it provides the only known channel to formulate a non-commutative Iwasawa main conjecture. At present, the only situation where the
$\M_H(G)$-conjecture is known to be valid is the ``$\mu=0$" situation (for instance,
see \cite[Proposition 5.6]{CF+} or \cite[Theorem 2.1]{CS12}). The
verification of the $\M_H(G)$-conjecture in general is still open (but see \cite[Section 2]{CSS}, \cite[Section 3]{CS12} or \cite[Section 3]{LimMHG} for
some related discussion in this direction; also see \cite{LimCom, LimMHGcong}).

For our purposes, we require the following.

\bl \label{MHG cyclotomic torsion}
Let $E$ be an elliptic curve defined over $F$ which has either good ordinary reduction or multiplicative reduction at each prime above $p$. Suppose that $F_\infty$ is a strongly admissible $p$-adic Lie extension of $F$ such that $X(E/F_\infty)$ belongs to $\M_H(G)$. Then for every finite extension $L$ of $F$ contained in $F_\infty$, the module $X(E/L^\cyc)$ is torsion over $\Zp\ps{\Gal(L^\cyc/L)}$.
\el

\bpf
 See \cite[Proposition 2.5]{CS12}.
\epf

\subsection{Akashi series of Selmer groups}

In this subsection, we review the calculation of Akashi series of Selmer group of an elliptic curve in a pro-$p$ extension. Such a calculation was first performed in \cite{CSS} and subsequently in \cite{CF+, Ze11}. The main result of this subsection is as follow, where we note that the calculations here does not require $\mathbf{(S2)}$.

\bp \label{akashi E formula}
Let $E$ be an elliptic curve defined over a number field $L$ which satisfies $\mathbf{(S1)}$. Suppose that $L_\infty$ is a strongly admissible pro-$p$ Lie extension of $L$ such that $X(E/L_\infty)\in\M_H(G)$, where $H=\Gal(L_\infty/L^\cyc)$ and $G=\Gal(L_\infty/L)$. Then we have
\begin{equation*}
\begin{split}
\Ak_H\big(X(E/L_\infty)\big)= & ~\ch_{\Zp\ps{\Ga}}\big(X(E/L^\cyc)\big)\times \prod_{w\in M(L^\cyc)}\prod_{j\geq 1}\ch_{\Zp\ps{\Ga}}\big(H^j(H_w, \Qp/\Zp)\big)^{(-1)^{j+1}} \\
 & \quad \times \prod_{w\in S'(L^\cyc)}\ch_{\Zp\ps{\Ga}}\big(H^1(H_w, \Ep(L_{\infty,w}))\big), \\
\end{split}
\end{equation*}
where $M(L^\cyc)$ denotes the set of primes of $L^\cyc$ above $p$ at which $E$ has split multiplicative reduction, and $S'(L^\cyc)$ is the set of primes of $L^\cyc$ above $S$ but not dividing $p$.
\ep

Before giving the proof of the proposition, we like to make the following remark.

\br
As noted above, such Akashi calculations have been performed in \cite{CF+, CSS, Ze11}. Here we should mention that our approach differs from these previous works. In these previous works, they have utilized a so-called ``large Selmer group" (see \cite[Formula (34)]{CSS} or \cite[Definition 3.1]{Ze11} for the precise definition of this Selmer group). Our adopted approach do not make use of this and is more direct. Furthermore, our approach has the advantage of removing the ``strongly admissible" assumption in \cite[Definition 1.2, Theorem 1.3]{Ze11}.
We also mentioned that we have give a detailed proof here, as some of the intermediate argument in the proof will be required for subsequent discussion of the main results in Sections \ref{main results} and \ref{twisted Section}.
\er

\bpf[Proof of Proposition \ref{akashi E formula}]
Plainly, the module $X(E/L_\infty)$ is torsion over $\Zp\ps{G}$. Furthermore, Lemma \ref{MHG cyclotomic torsion} tells us that $X(E/L^\cyc)$ is torsion over $\Zp\ps{\Ga}$. In view of these observations, it follows from Propositions \ref{cyclotomic torsion} and \ref{p-adic torsion} that we have a short exact sequence
\[ 0\lra  \Sel(E/\mathcal{L})\lra H^1(G_S(\mathcal{L}),E[p^\infty])\lra \bigoplus_{v\in S}J_v(E/\mathcal{L})\lra 0\]
for $\mathcal{L}= L^\cyc, L_\infty$. The short exact sequence for $L^\cyc$ and the $H$-cohomology long exact sequence associated to the corresponding short exact sequence for $L_\infty$ fit into the following commutative diagram
 \[   \entrymodifiers={!! <0pt, .8ex>+} \SelectTips{eu}{}\xymatrix{
    0 \ar[r]^{} &\Sel(E/L^\cyc)  \ar[d]^{\alpha}  \ar[r] &  H^1\big(G_{S}(L^\cyc),\Ep\big)
    \ar[d]^{\beta} \ar[r] & \displaystyle \bigoplus_{v\in S}J_v(E/L^\cyc) \ar[d]^{r=\oplus r_v} \ar[r]& 0 \\
    0 \ar[r]^{} & \Sel(E/L_{\infty})^{H} \ar[r]^{} & H^1\big(G_{S}(L_{\infty}),\Ep\big)^{H} \ar[r] &
    \left(\displaystyle \bigoplus_{v\in S}J_v(E/L_{\infty})\right)^{H} \ar[r] &~ \cdots  } \]
    with exact rows. From which, we obtain a long exact sequence
\begin{equation*}
\begin{split}
0\lra \ker\al \lra \ker\be \lra \ker r \lra  \coker \al  \lra \coker \be \lra \coker r \lra H^1(H, \Sel(E/L_\infty))\\
 \lra H^1\left(H,H^1\big(G_{S}(L_{\infty}),\Ep\big)\right) \lra H^1\left(H,\bigoplus_{v\in S}J_v(E/L_{\infty})\right) \lra \cdots \\
\end{split}
\end{equation*}

By Propositions \ref{cyclotomic torsion} and \ref{p-adic torsion} again, one has $H^2(G_S(L^\cyc),\Ep)=H^2(G_S(L_\infty),\Ep)=0$. Therefore, the spectral sequence
\[H^i\big(H, H^j(G_S(L_\infty),\Ep)\big)\Longrightarrow H^{i+j}\big(G_S(L^\cyc),\Ep\big)\]
degenerates to yield an exact sequence
\[ 0 \lra H^1\big(H, \Ep(L_\infty)\big) \lra H^j\big(G_S(L^\cyc),\Ep\big) \lra H^j\big(G_S(L_\infty),\Ep\big)^H \]\[\lra H^2\big(H, \Ep(L_\infty)\big) \lra 0\]
and isomorphisms
\[ H^j\big(H, H^1(G_S(L_\infty),\Ep)\big) \cong H^{j+2}\big(H, \Ep(L_\infty)\big)\]
for $j\geq 1$. In particular, this shows that $\ker \be$, $\coker \be$ and $H^j\big(H, H^1(G_S(L_\infty),\Ep)\big)$ are cofinitely generated over $\Zp$. Furthermore, their characteristic elements (over $\Zp\ps{\Ga}$) can be computed in terms of the characteristic elements of $H^j\big(H, \Ep(L_\infty)\big)$.

Let $v$ be a prime of $F$ which lies above $p$. By Lemma \ref{local coh description}, we have $J_v(E/L^{\cyc}) = \bigoplus_{w|v} H^1(L^\cyc_w, D_w)$, where the sum is over the primes of $L^\cyc$ above $w$. We then write $H_w$ for the decomposition group of $H$ at some fixed prime of $L_\infty$ above $w$. It follows from an application of Shapiro's lemma that
\[H^j\big(H, J_v(E/L_\infty)\big)\cong \bigoplus_{w|v}H^j\big(H_w, H^1(L^\cyc_w, D_w)\big).\]
Since $H^2\big(L^\cyc_w, D_w\big)=H^2\big(L_{\infty,w}, D_w\big)=0$ by \cite[Theorem 7.1.8(i)]{NSW}, we may apply a similar proof as above to obtain an exact sequence
\[ 0 \lra H^1(H_w, D_w(L_{\infty,w})) \lra H^1(L^\cyc_w, D_w) \lra H^1(L_{\infty,w}, D_w)^{H_w} \]\[\lra H^2(H_w, D_w(L_{\infty,w})) \lra 0\]
and isomorphisms
\[ H^j\big(H_w, H^1(L_{\infty,w}, D_w)\big) \cong H^{j+2}\big(H_w, D_w(L_{\infty,w})\big)\]
for $j\geq 1$. One can perform similar calculations for the local terms at primes outside $p$.
Since $\Zp\ps{\Ga}$-characteristic elements are multiplicative in short exact sequences of torsion $\Zp\ps{\Ga}$-modules, and taking Lemma \ref{akashi ord p} into account, we have
\begin{equation*}
\begin{split}
\Ak_H(X(E/L_\infty)= & ~\ch_{\Zp\ps{\Ga}}\big(X(E/L^\cyc)\big)\times \prod_{w\in M(L^\cyc)}\prod_{j\geq 1}\ch_{\Zp\ps{\Ga}}\big(H^j(H_w, \Qp/\Zp)\big)^{(-1)^{j+1}} \\
 & \quad \times \prod_{w\in S'(L^\cyc)}\ch_{\Zp\ps{\Ga}}\big(H^1(H_w, \Ep(L_{\infty,w}))\big)\times \prod_{j\geq 1}\ch_{\Zp\ps{\Ga}}\left( H^j(H, \Ep(L_\infty))\right)^{(-1)^j} \\
\end{split}
\end{equation*}
Here we also note that in view of Lemma \ref{Non-split local}, the local terms at non-split multiplicative primes above $p$ have no contributions to $\Zp\ps{\Ga}$-characteristic elements, and so they do not appear in the above formula.
The conclusion of the proposition is now a consequence of this and the next lemma.
\epf

\bl \label{akashi E global}
 Let $E$ be an elliptic curve defined over a number field $L$ which satisfies $\mathbf{(S1)}$. Let $L_\infty$ be a strongly admissible pro-$p$ $p$-adic Lie extension of $L$.
Then we have the following equalities.
 \begin{enumerate}
 \item[$(i)$] $\Ak_H(\Ep(L_\infty))=1$.
 \item[$(ii)$] $\displaystyle\prod_{j\geq 1}\ch_{\Zp\ps{\Ga}}\left( H^j(H, \Ep(L_\infty))\right)^{(-1)^j} =1$.
 \end{enumerate}
\el

\bpf
By \cite[Lemma 5.3]{Ze09}, if $\Ep$ is not realized over $L_\infty$, then $\Ep(L_\infty)$ is finite. Therefore, the assertion of the lemma is clear in this situation. Now suppose that $\Ep(L_\infty) = \Ep$, then $L_\infty$ contains $L(\Ep)$. Since $L_\infty/L$ is pro-$p$, so is $L(\Ep)/L$. Hence we have \[\Gal(L(\Ep)/L)\cong \Gal(L(\Ep)/L^\cyc) \times \Gal(L^\cyc/L),\] where $\Gal(L(\Ep)/L^\cyc)$ is isomorphic to $\Zp$ or an open pro-$p$ subgroup of $\mathrm{SL}_2(\Zp)$ accordingly to $E$ having complex multiplication or not. Either way, Lemma \ref{Akashi Zp} applies.
\epf

\section{First main result} \label{main results}

In this section, we will present and prove our first main result. As a start, we introduce one more hypothesis.

\begin{itemize}
\item[$\mathbf{(S3)}$] For each $v$ of $F$ above $p$ at which $E$ has non-split multiplicative reduction, assume that for every finite extension $L$ of $F$ contained in $F_\infty$, the elliptic curve $E$ has non-split multiplicative reduction at every prime of $L$ above $v$.
\end{itemize}

We introduce one last hypothesis to handle the primes of split multiplication reduction of $E$ above $p$.

\bd \label{splitmultdef}
Let $E$ be an elliptic curve defined over a number field $F$ and $F_\infty$ a strongly admissible $p$-adic Lie extension of $F$. An extension $L$ of $F$ contained in $F_\infty$ is said to satisfy $\mathbf{(M_p)}$ if for each prime $w$ of $L$ above $p$ at which $E$ has split multiplicative reduction, either of the following holds.

\begin{itemize}
     \item[$(\mathrm{I})$] For every prime $x$ of $F_\infty$ above $w$, the extension $F_{\infty,x}$ contains a subextension $\mathcal{L}_\infty$ of $L_w^{\cyc}$ such that  $\Gal(\mathcal{L}_\infty/L_w)\cong \Gal(\mathcal{L}_\infty/L_w^\cyc)\times \Gal(L_w^\cyc/L_w)$, where $\Gal(\mathcal{L}_\infty/L_w^\cyc)$ has dimension $\geq 1$.

     \item[$(\mathrm{II})$] $L_w$ contains a primitive $p$-root of unity and $L_{\infty, x}$ is a multi-False-Tate extension over $L_w$.
\end{itemize}

We shall write $m_p(L)$ for the number of primes of $L^\cyc$ (above $p$ and at which $E$ has split multiplicative reduction) satisfying statement (I).
\ed

\br
Hypothesis (S3) is more of a convenience to simplify the presentation. However, the assumptions in Definition \ref{splitmultdef} are crucial, as these constraints are put upon us by Lemmas \ref{akashi split mult 1} and \ref{akashi zero order split mult 2}.
\er

\bt \label{main: pro-p}
Let $E$ be an elliptic curve defined over a number field $F$ which satisfies $\mathbf{(S1)-(S3)}$, and let $F_\infty$ be a strongly admissible $p$-adic Lie extension of $F$ such that $X(E/F_\infty)\in \M_H(G)$. Suppose that $L$ is a finite Galois extension of $F$ contained in $F_\infty$ which satisfies the following three statements.
\begin{enumerate}
\item[$(a)$] $F_\infty/L$ is a pro-$p$ extension.
\item[$(b)$] Conjecture \ref{semisimple conj} is valid for $X(E/L^\cyc)$.
\item[$(c)$] $L$ satisfies $\mathbf{(M_p)}$.
\end{enumerate}
If $\xi_E$ is a characteristic element of $X(E/F_\infty)$ in the sense of Definition \ref{CoatesFKSV}, we then have
\[\displaystyle\ord_{T=0} ~\big(\Phi_{\reg_L}(\xi_E)\big) = \corank_{\Zp}\big(\Sel(E/L)\big) + m_p(L),\]
where $m_p(L)$ is defined as in Definition \ref{splitmultdef}.
\et

\bpf
By Lemma \ref{independence of generator and units}, the term on the left is independent of the choice of $\xi_E$. Therefore, we may simply work with one such fixed choice. Write $U=\Gal(F_\infty/L)$ and $\Ga_L= \Gal(L^\cyc/L)$. Fix an appropriate power of the generator of $\Ga$ such that it is a generator for $\Ga_L$. Under these choices of generators, we have the following identifications and inclusion
\[ \Zp\ps{T_L}\cong \Zp\ps{\Ga_L} \subseteq \Zp\ps{\Ga}\cong \Zp\ps{T},\]
 where $T_L = (1+T)^{|\Ga:\Ga_L|} -1$. Recall from Subsection \ref{char element subsec} that we have a restriction map on the $K_1$-groups
\[\res : K_1(\Zp\ps{G}_{\Si^*}) \lra K_1(\Zp\ps{U}_{\Si^*_U}). \]
By Lemma \ref{res BVseries}, we see that $\res(\xi_E)$ is a characteristic element of $X(E/F_\infty)$ viewed in $\M_{H_L}(U)$, where $H_L = \Gal(F_\infty/L^\cyc)$. Taking this into account, it then follows from Proposition \ref{evaluation ind rep} and Lemma \ref{char and Akashi} that $\Phi_{\reg_L}(\xi_E) = \Ak_{H_L}\big(X(E/F_\infty)\big)$ mod $\Op\ps{T_L}^{\times}$. By Proposition \ref{akashi E formula}, the latter is given by
\begin{equation*}
\begin{split}
\Ak_{H_L}\big(X(E/F_\infty)\big)= & ~\ch_{\Zp\ps{\Ga_L}}\big(X(E/L^\cyc)\big)\times \prod_{w\in M(L^\cyc)}\prod_{j\geq 1}\ch_{\Zp\ps{\Ga_L}}\big(H^j(H_w, \Qp/\Zp)\big)^{(-1)^{j+1}} \\
 & \quad \times \prod_{w\in S'(L^\cyc)}\ch_{\Zp\ps{\Ga_L}}\big(H^1(H_w, \Ep(L_{\infty,w}))\big), \\
\end{split}
\end{equation*}
where $M(L^\cyc)$ is the set of primes of $L^\cyc$ above $p$ at which $E$ has split multiplicative reduction. It remains to calculate the order of vanishing of the terms appearing on the right.  Lemma \ref{neq l structure twist} tells us that the local cohomology terms outside $p$ has no contribution to order of vanishing. For the split multiplicative primes above $p$ and in view of assumption (c), we may invoke Lemmas \ref{akashi split mult 1} and \ref{akashi zero order split mult 2} to calculate their contributions to the order of vanishing. In particular, each $w\in M(F^\cyc)$ contribute either 1 or 0 accordingly to whether the prime $w$ satisfies (I) or (II) in Definition \ref{splitmultdef}. Finally, by Proposition \ref{semisimple} and assumption (b), the order of vanishing of $\ch_{\Zp\ps{\Ga_L}}\big(X(E/L^\cyc)\big)$ is given by $\corank_{\Zp}\big(\Sel(E/L)\big)$. Combining these observations, we have our theorem.
\epf

\subsection{Appendix: Relation with generalized Euler characteristics}

In this short appendix, we describe how the characteristic element is related to the cyclotomic generalized Euler characteristics.  This relation is somewhat documented in \cite{CSS, Ze09}. We will mainly be concerned of the contribution of our Theorem \ref{main: pro-p} towards this. For simplicity, we shall assume that our elliptic curve $E$ has good ordinary reduction at all primes above $p$. For a discussion when $E$ has multiplicative reduction, we refer readers to \cite{DL}. As a start, we recall the following theorem obtained by Perrin-Riou \cite{PR} and Schneider \cite{Sch85} independently.

\bt[Perrin-Riou, Schneider]
Let $E$ be an elliptic curve defined over a number field $L$ which has good ordinary reduction at all primes above $p$. Suppose that $X(E/L^\cyc)$ is torsion. Furthermore, assume that $\sha(E/L)[p^\infty]$ is finite. Then the leading coefficient of $\ch_T\big(X(E/L^\cyc)\big)$ is given by
\[ \mathrm{Reg}_p(E/L) \times
\frac{|\sha(E/L)_{p^\infty}|}{|\Ep(L)|^2} \times  \prod_v c_v^{(p)} \times \prod_{v\in S_p}(d_v^{(p)})^2.\]
 Here $\mathrm{Reg}_p(E/L)$ is the normalised $p$-adic regulator on $E(L)$, $c_v^{(p)}$ is the highest power of $p$ dividing $|E(L_v):E_0(L_v)|$, where $E_0(L_v)$ is the subgroup of $E(L_v)$ consisting of points with nonsingular reduction modulo $v$, and $d_v^{(p)}$ is the highest power of $p$ dividing $|\tilde{E}_v(l_v)|$, where $l_v$ is the residue field of $L_v$.
\et

In the event that the elliptic curve $E$ satisfies Conjecture \ref{semisimple conj}, it then follows from Proposition \ref{semisimple} that the leading coefficient of $\ch_T\big(X(E/L^\cyc)\big)$ can be expressed as
\[ \frac{1}{T^{r(L)}} ~\ch_T\big(X(E/L^\cyc)\big)\Big|_{T=0},  \]
where $r(L) = \corank_{\Zp}\big(\Sel(E/L)\big)$.
Now combining the above theorem with Propositions \ref{evaluation ind rep} and \ref{akashi E formula}, we obtain the following.

\bp \label{main: pro-p relate euler}
Let $E$ be an elliptic curve defined over a number field $F$  which has good ordinary reduction at all primes above $p$, and let $F_\infty$ be a strongly admissible $p$-adic Lie extension of $F$ such that $X(E/F_\infty)\in \M_H(G)$. Suppose that $L$ is a finite Galois extension of $F$ contained in $F_\infty$ which satisfies the following three statements.
\begin{enumerate}
\item[$(a)$] $F_\infty/L$ is a pro-$p$ extension.
\item[$(b)$] Conjecture \ref{semisimple conj} is valid for $X(E/L^\cyc)$.
\item[$(c)$] $\sha(E/L)_{p^\infty}$ is finite.
\end{enumerate}
If $\xi_E$ is a characteristic element of $X(E/F_\infty)$ and $r(L) = \rank_{\Zp}\big(X(E/L)\big)$, we then have
\[\frac{1}{T^{r(L)}} ~\big(\Phi_{\reg_L}(\xi_E)\big)\Big|_{T=0} \hspace{5in}\]
\[= |\Ga:\Ga_L|^{r(L)}\times\mathrm{Reg}_p(E/L) \times\frac{|\sha(E/L)_{p^\infty}|}{|\Ep(L)|^2} \times  \prod_v c_v^{(p)} \times \prod_{w\in S_p(L)}(d_v^{(p)})^2 \times \left|\prod_{w\in \mathcal{R}_L}L_w(E,1)\right|_p ,\]
where $\mathcal{R}_L$ is the set of primes of $L$ which do not divide $p$ but ramify in $F_\infty/L$.
\ep

\br
As mentioned in the beginning of this appendix, the formula is somewhat proven in \cite{CSS, Ze09}. The contribution of our Theorem \ref{main: pro-p} comes in two forms. One lies in the order of vanishing, namely, the term $\frac{1}{T^{r(L)}}$ on the left of the above equation. The second contribution is the term $|\Ga:\Ga_L|^{r(L)}$ on the right.
\er

\section{Artin twist of characteristic element} \label{twisted Section}

In this section, we consider some cases where the order of vanishing of the characteristic element at Artin twist can be evaluated.

\subsection{Twist of Selmer groups} \label{twisted Selmer Subsec}

Retain the settings in Subsection \ref{p-adic Lie subsection}. Let $\rho:G\lra \mathrm{GL}_d(\Op)$ be an Artin representation of $G$ and write $W_\rho$ for a free $\Op$-module of rank $d$ realizing $\rho$. For any extension $\mathcal{L}$ of $F^\cyc$ contained in $F_\infty$, the twisted Selmer group of $E$ over $\mathcal{L}$ (\cite{CFKS, G89}) is defined by
\[ \Sel(\tw_{\rho}(E)/\mathcal{L})=\ker\Big(H^1(G_S(\mathcal{L}),W_\rho\ot_{\Zp}\Ep)\lra \bigoplus_{v\in S}J_v(\tw_{\rho}(E)/\mathcal{L})\Big),\]
where
\[ J_v(\tw_{\rho}(E)/\mathcal{L})\cong \begin{cases} \ilim_{\mathcal{L}'}\bigoplus_{w|v}H^1(\mathcal{L}'_w, W_\rho\ot_{\Zp}D_w),  & \mbox{if $v$ divides $p$}, \\
\ilim_{\mathcal{L}'}\bigoplus_{w|v}H^1(\mathcal{L}'_w, W_\rho\ot_{\Zp}\Ep), & \mbox{if $v$ does not divide $p$}.\end{cases} \]

The Pontryagin dual of $ \Sel(\tw_{\rho}(E)/\mathcal{L})$ is denoted by $X(\tw_{\rho}(E)/\mathcal{L})$. The arithmetic significance of these twisted Selmer groups lies in the following definition and lemma.

\bd[\cite{CFKS, Lee}]
Suppose that $\rho$ is irreducible. Write $\mF$ for any finite Galois extension of $F$ contained in $F_\infty$ such that $\rho$ factors through $\Gal(\mF/F)$. We then define $s_{E,\rho}$ to be the number of copies of $W_\rho\ot_{\Op}\bar{\Q}_p$ occurring in $X(E/\mF)\ot_{\Zp}\bar{\Q}_p$.
\ed

\bl \label{CFKSroot}
If $\rho$ is an irreducible representation of $G$, the following statements are valid.
\begin{enumerate}
\item[$(i)$] $W_{\rho}\ot_{\Zp}X(E/F_\infty)  = X\big(\tw_{\hat{\rho}}(E)/F_\infty\big)$, where $\hat{\rho}$ is the contragradient of $\rho$.
\item[$(ii)$] $s_{E,\rho} = \corank_{\Op}\big(\Sel\big(\tw_{\rho}(E)/F\big)\big)$.
\end{enumerate}
\el

\bpf
See \cite[Lemma 3.4]{CFKS}.
\epf

The next lemma records some consequences of $X(E/F_\infty) \in \M_H(G)$ on the twisted Selmer groups.

\bl \label{twisted tor property}
 Suppose that $X(E/F_\infty) \in \M_H(G)$. Then $X(\tw_{\rho}(E)/F_\infty) \in \M_H(G)$ and for every finite extension $L$ of $F$ contained in $F_\infty$, we have $X(\tw_{\rho}(E)/L^\cyc)$ is torsion over $\Op\ps{\Ga_L}$, where $\Ga_L=\Gal(L^\cyc/L)$. Furthermore, we have short exact sequences
 \[ 0\lra \Sel(\tw_{\rho}(E)/\mathcal{L}) \lra H^1(G_S(\mathcal{L}),W_\rho\ot_{\Zp}\Ep)\lra \bigoplus_{v\in S}J_v(\tw_{\rho}(E)/\mathcal{L}) \lra 0\]
 and $H^2(G_S(\mathcal{L}),W_\rho\ot_{\Zp}\Ep)=0$ for $\mathcal{L} = F_\infty, L^\cyc$.
\el

\bpf
  This lemma and its proof should be known among experts. Due to a lack of proper reference, we shall give a sketch of the ideas behind it. As a start, one recalls the well-known fact that the property of a $\Zp\ps{G}$-module lying in $\M_H(G)$ is invariant under Artin twists (cf. \cite[Lemma 3.2]{CF+}). Combining this fact with Lemma \ref{CFKSroot} and the hypothesis that $X(E/F_\infty) \in \M_H(G)$, we have the assertion that $X(\tw_{\rho}(E)/F_\infty) \in \M_H(G)$. In view of this $\M_H(G)$-property, we may then apply a descent argument similar to that in \cite[Proposition 2.5]{CS12} to obtain
  the $\Op\ps{\Ga_L}$-torsionness of $X(\tw_{\rho}(E)/L^\cyc)$. Taking this latter torsionness observation into account, one can invoke a standard argument, which consists of global/local Euler characteristics calculations and an analysis of the Poitou-Tate sequence (for instances, see \cite[Proposition 3.4 and Corollary 3.5]{LimMHG} or \cite[Lemma 5.1.2]{LimMHGcong}), to obtain the remaining assertions of the lemma.
\epf

For subsequent discussion, we shall fix a uniformizer $\pi$ for $\Op$. We now state the following analogue of Greenberg's semisimple conjecture for our twisted Selmer groups.

\begin{conjecture}[Greenberg] \label{semisimple conj2}
Under torsion hypothesis, $X(\tw_{\rho}(E)/F^\cyc)$ satisfies Greenberg's semisimple conjecture. In other words, we have a pseudo-isomorphism
\[X(\tw_{\rho}(E)/F^\cyc) \sim \bigoplus_{i=1}^s\Op\ps{\Ga}/\pi^{\al_i} \times \bigoplus_{j=1}^t\Op\ps{\Ga}/f_j,\]
where each $f_j$ is irreducible in $\Op\ps{\Ga}$ and is not an associate of $\pi$.
\end{conjecture}

The following theorem describes the order of vanishing of the characteristics element evaluated at Artin representations as introduced in Theorem \ref{main: Artin twist intro}.

\bt \label{main: Artin twist}
Let $E$ be an elliptic curve defined over a number field $F$ which has good ordinary reduction at all primes of $F$ above $p$. Let $F_\infty$ be a admissible $p$-adic Lie extension of $F$. Suppose that all of the following statements are valid.
\begin{enumerate}
\item[$(a)$] $X(E/F_\infty)\in \M_H(G)$.
\item[$(b)$] Conjecture \ref{semisimple conj2} is valid for $ X(\tw_{\rho}(E)/F^\cyc)$.
\item[$(c)$] For every open subgroup $H'$ of $H$, $H^i(H', \Ep(F_\infty))$ is finite for all $i\geq 1$.
 \item[$(d)$] For every prime $w$ of $F^\cyc$ dividing $p$, and each open subgroup $H_w'$ of $H_w$, $H^i(H_w', \widetilde{E}_{p^\infty}(F_\infty))$ is finite for all $i\geq 1$.
\end{enumerate}
Let $\xi_E$ be a characteristic element of $X(E/F_\infty)$ and $\rho$ an irreducible Artin representation of $G=\Gal(F_\infty/F)$. Then
we have
\[\displaystyle\ord_{T=0} ~\big(\Phi_{\rho}(\xi_E)\big) = s_{E,\rho}.\]
\et

\bpf The proof of this follows closely to that in Theorem \ref{main: pro-p}. In view of Lemma \ref{twisted tor property}, we have the following commutative diagram
 \[   \entrymodifiers={!! <0pt, .8ex>+} \SelectTips{eu}{}\xymatrix{
    0 \ar[r]^{} &\Sel(\tw_\rho(E)/F^\cyc)  \ar[d]^{\alpha}  \ar[r] &  H^1\big(G_{S}(F^\cyc),W_\rho\ot_{\Op}\Ep\big)
    \ar[d]^{\beta} \ar[r] & \displaystyle \bigoplus_{v\in S}J_v(\tw_\rho(E)/F^\cyc) \ar[d]^{r=\oplus r_v} \ar[r]& 0 \\
    0 \ar[r]^{} & \Sel(\tw_\rho(E)/F_{\infty})^{H} \ar[r]^{} & H^1\big(G_{S}(F_{\infty}),W_\rho\ot_{\Op}\Ep\big)^{H} \ar[r] &
    \left(\displaystyle \bigoplus_{v\in S}J_v(\tw_\rho(E)/F_{\infty})\right)^{H} \ar[r] &~ \cdots  } \]
    with exact rows. By a snake lemma argument, we obtain a long exact sequence
\begin{equation*}
\begin{split}
0\lra \ker\al \lra \ker\be \lra \ker r \lra  \coker \al  \lra \coker \be \lra \coker r \lra H^1(H, \Sel(E/F_\infty))\\
 \lra H^1(H,W_\infty) \lra H^1(H,J_\infty) \lra H^2(H, \Sel(E/F_\infty))\lra \cdots \\
\end{split}
\end{equation*}
where $W_\infty =  H^1\big(G_{S}(F_{\infty}),W_{\rho}\ot_{\Zp}\Ep\big)$ and $J_\infty = \bigoplus_{v\in S}J_v(\tw_{\rho}(E)/F_{\infty})$. By a similar argument to that in Proposition \ref{akashi E formula} and taking Lemma \ref{twisted tor property} into account, we have that $\ker \be$, $\coker \be$ and $H^i\big(H, W_\infty\big)$ are cofinitely generated over $\Op$, and their characteristic elements (over $\Op\ps{\Ga}$) can be computed in terms of the characteristic elements of $H^i\big(H, W_{\rho}\ot\Ep(F_\infty)\big)$. Via a similar argument to that in Lemma \ref{akashi ord p twist}, the order of vanishing at $T=0$ for these latter elements is zero. Similarly, by appealing to Lemmas \ref{akashi ord p twist} and \ref{neq l structure twist}, the $\Op\ps{\Ga}$-characteristic elements of
$(\ker \ga)^\vee$, $(\coker \ga)^\vee$ and $H^i\big(H, J_\infty\big)^\vee$ (for $i\geq 1$) have trivial orders of vanishing at $T=0$. Consequently, the $\Op\ps{\Ga}$-characteristic elements of
$(\ker \al)^\vee$ and $(\coker \al)^\vee$, and $H^i\big(H, \Sel(\tw_{\rho}(E)/F_\infty)\big)^\vee$ (for $i\geq 1$) have trivial orders of vanishing at $T=0$. Combining these observations with Lemmas \ref{char and Akashi} and \ref{CFKSroot}(i), we have
\[\displaystyle\ord_{T=0}\big(\Phi_{\rho}(\xi_E)\big) = \ord_{T=0}\big(\Ak_H(\tw_{\hat{\rho}}X(E/F_\infty))\big) = \ord_{T=0}\big(\Ak_H(X(\tw_{\rho}(E)/F_\infty))\big) \]
\[=\ord_{T=0} \big(\ch_{\Op\ps{\Ga}}X(\tw_{\rho}(E)/F_\infty)_H\big) =\ord_{T=0} \big(\ch_{\Op\ps{\Ga}}X(\tw_{\rho}(E)/F^\cyc)\big). \]
Now, a similar argument to that in Proposition \ref{semisimple} shows that
\[ \Sel(\tw_{\rho}(E)/F) \lra \Sel(\tw_{\rho}(E)/F^\cyc)^{\Ga}\]
has finite kernel and cokernel. As in Proposition \ref{semisimple}, combining this latter fact with hypothesis (b), we have
\[\ord_{T=0} \big(\ch_{\Op\ps{\Ga}}X(\tw_{\rho}(E)/F^\cyc)\big) = \corank_{\Op}\big(\Sel\big(\tw_{\rho}(E)/F\big)\big), \]
where the latter is $s_{E,\rho}$ by Lemma \ref{CFKSroot}(ii). This proves the proposition.
\epf

\br \label{hypotheses cd}
The hypotheses (c) and (d) are rather mild. In fact, they are known to be satisfied for many extensions.

(i) When dim $G \leq 3$, hypotheses (c) and (d) are verified in the proof of \cite[Lemma 2.3]{DL15}.

(ii) Let $F_\infty = \Q\big(\mu_{p^\infty}, \sqrt[p^\infty]{\alpha_1} , \ldots, \sqrt[p^\infty]{\alpha_{d-1}}\big)$, where $\alpha_1, \ldots , \alpha_{d-1}\in \Q^{\times}$, whose image in $\Q^\times/(\Q^{\times})^p$ are linearly independent over $\Z/p\Z$. In this situation, the validity of the hypotheses follow from the result of Kubo-Taguchi \cite[Theorem 1.1]{KT}.

(iii) If $F_\infty = F(E_{p^\infty})$, where $E$ is the elliptic curve in question, the validity of the hypotheses has been addressed in \cite[Corollary 1.4]{CSW}.

(iv) In \cite[Section 5]{Ze09}, the hypotheses have also been verified for a large class
of $p$-adic Lie extensions.
\er

\subsection{False-Tate extension: an alternative approach} \label{false-tate subsec}

In this subsection, we consider the case of a False-Tate extension. As before, $p$ denotes an odd prime. Let $m$ be a $p$-powerfree integer. Set $F_\infty = \Q(\mu_{p^\infty}, \sqrt[p^\infty]{m})$. This is a Galois extension of $\Q$ with Galois group $G=\Gal(F_\infty/\Q)=\Zp\rtimes \Zp^{\times}$. Write $H=\Gal(F_\infty/\Q^\cyc)$ and $\Ga= \Gal(\Q^\cyc/\Q)$. For $n\geq 1$, let $\rho_n$ denote the representation of $G$ obtained by inducing
any character of exact order $p^n$ of $\Gal\big(\Q(\mu_{p^n}, \sqrt[p^n]{m})/\Q(\mu_{p^n})\big)$ to $\Gal\big(\Q(\mu_{p^n}, \sqrt[p^n]{m})/\Q\big)$.
Let $E$ be an elliptic curve defined over $\Q$. We can now prove the following final result as stated in the introduction.

\bt \label{ArtinFT}
Let $E$ be an elliptic curve defined over a number field $\Q$ which has either good ordinary reduction or multiplicative reduction at $p$. Suppose that $X(E/F_\infty)\in \M_H(G)$. Assume that Conjecture \ref{semisimple conj} is valid for $X(E/L_n^\cyc)$, where $L_n = ~\Q(\mu_{p^n}, \sqrt[p^n]{m})$.

If $\xi_E$ is a characteristic element of $X(E/F_\infty)$, we then have
\[\displaystyle\ord_{T=0} ~\big(\Phi_{\rho_n}(\xi_E)\big) = s_{E,\rho_n}.\]
\et

\bpf
To lighten notation, we shall write $\rho = \rho_n$. Since $\rho$ is self-dual, Lemma \ref{char and Akashi} reads as
 \[ \Phi_{\rho}(\xi_E) = \Ak_H\big(\tw_{\rho}(X(E/F_\infty))\big)~(\mbox{mod} ~\Zp\ps{\Ga}^\times). \]
 To continue, we introduce more notations. Set $L_n'=\Q(\mu_{p^n}, \sqrt[p^{n-1}]{m})$. Let $G_n = \Gal(F_\infty/L_n)$ and $G'_n = \Gal(F_\infty/L_n')$. Write $H_n= G_n\cap H$ and $H'_n=G'_n\cap H$. Note that $G_n/H_n = G'_n/H'_n$ identifies with $\Gal(\Q(\mu_{p^n})/\Q(\mu_p))$ which we shall denote by $\Ga'$. If we identify $\Ga=\Gal(\Q^\cyc/\Q)$ with $\Gal\big(\Q(\mu_{p^\infty})/\Q(\mu_p)\big)$, then $\Ga'$ can be viewed as a subgroup of $\Ga$ of index $p^{n-1}$. With these notations in hand, Proposition \ref{evaluation ind rep} tells us that $\Phi_{\rho}(\xi_E)$ lies in $Q_{\Op}(\Ga')$. Now, by virtue of Lemma \ref{char and Akashi}, we may write
 \[ \Ak_H\big(\tw_{\rho}(X(E/F_\infty))\big) = \Phi_{\rho}(\xi_E)\cdot u\]
 for some $u \in \Zp\ps{\Ga}^{\times}$.
On the other hand, it follows from \cite[Lemma A.18]{DD} that we have
\[ \Ak_{H'_n}\big(X(E/F_\infty)\big) = \Ak_{H'_n}\big(X(E/F_\infty)\big)N_{\Ga/\Ga'}\big(\Ak_H(\tw_{\rho}(X(E/F_\infty)))\big)^{p-1},  \]
where $N_{\Ga/\Ga'}$ is the norm map from $\Zp\ps{\Ga}$ to $\Zp\ps{\Ga'}$. Combining these observations, we obtain
\begin{equation}\label{Artin formalism for Ak}
  \Ak_{H'_n}\big(X(E/F_\infty)\big) = \Ak_{H'_n}\big(X(E/F_\infty)\big)\Phi_{\rho}(\xi_E)^{p^{n-1}(p-1)}
N_{\Ga/\Ga'}(u)^{p-1}.
\end{equation}
By the hypothesis of the proposition, we may apply Proposition \ref{char and Akashi} and Theorem \ref{main: pro-p} to see that
\[ \ord_{T=0}\Big(\Ak_{H_n}\big(X(E/F_\infty)\big)\Big) = \corank_{\Zp} \big(\Sel(E/L_n)\big). \]
(Note that we are in the situation (II) of Definition \ref{splitmultdef}.) We now claim that $X(E/L_n'^{\cyc})$ satisfies Conjecture \ref{semisimple conj} as a $\Zp\ps{\Ga'}$-module. Indeed, via a descent argument, one can show that the map
\[ X(E/L_n^\cyc)_{\Gal(L_n^\cyc/L_n'^\cyc)} \lra X(E/L_n'^{\cyc})\]
has finite kernel and cokernel. Since Conjecture \ref{semisimple conj} is valid for $X(E/L_n^\cyc)$, it follows from the above that the same can be said for $X(E/L_n'^{\cyc})$. This proves our claim. Consequently, we may apply Proposition \ref{char and Akashi} and Theorem \ref{main: pro-p} to conclude that
\[ \ord_{T=0}\Big(\Ak_{H'_n}\big(X(E/F_\infty)\big)\Big) = \corank_{\Zp} \big(\Sel(E/L_n')\big). \]
Putting these into (\ref{Artin formalism for Ak}), we obtain
\begin{equation}\label{Phi rank}
   p^{n-1}(p-1)\ord_{T=0} ~\big(\Phi_{\rho_n}(\xi_E)\big) =  \corank_{\Zp} \big(\Sel(E/L_n)\big) - \corank_{\Zp} \big(\Sel(E/L_n')\big).
\end{equation}
On the other hand, the set of all irreducible representations of $\Gal(L_n/\Q)$ consists of $\rho$ together with all
irreducible representations of $\Gal(L_n'/\Q)$. Consequently, it follows that
\[ X(E/L_n)\ot\bar{\Q}_p = (X(E/L_n')\ot\bar{\Q}_p) \oplus (W_\rho\ot\bar{\Q}_p)^{\oplus s_{E,\rho}}  \]
which in turn implies that
\begin{equation}\label{s rank}
\corank_{\Zp} \big(\Sel(E/L_n)\big) - \corank_{\Zp} \big(\Sel(E/L_n')\big) = s_{E,\rho} \rank_{\Zp} (W_\rho) =  p^{n-1}(p-1)s_{E,\rho}.
\end{equation}
Comparing the two equalities (\ref{Phi rank}) and  (\ref{s rank}), we have the conclusion of the theorem.
\epf

\section{Examples} \label{examples section}

We now give some classes of examples to illustrate our results.

\subsection{Good ordinary consideration}

To facilitate our discussion, we need to recall the following result (cf. \cite[Proposition A.37]{DD}).

\bp \label{DD CS}
Let $E$ be an elliptic curve defined over $\Q$ with good ordinary reduction at $p$, and let $q$ be a multiplicative reduction prime of $E$. Suppose that one of the following statements holds.
\begin{enumerate}
  \item[$(a)$] $E$ has split multiplicative reduction at the prime $q$ with $q$ being inert in $\Q(\mu_p)/\Q$ and $\Sel(E/\Q(\mu_{p^\infty})) = 0$.
  \item[$(b)$] $E$ has non-split multiplicative reduction at the primes of $\Q(\mu_p)$ above $q$ and $\Sel(E/\Q(\mu_{p^\infty})) = \Qp/\Zp$ with a trivial $\Ga$-action.
\end{enumerate}
 Then $X(E/F_\infty)$ is a free $\Zp\ps{\Gal(F_\infty/\Q(\mu_{p^\infty}))}$-module of rank 1, where $F_\infty = \Q(\mu_{p^{\infty}},\sqrt[p^\infty]{q})$. \ep

We now establish the following.

\bp \label{uncond}
Retain settings in Proposition \ref{DD CS}. Then $X(E/F_\infty)\in\M_H(G)$. Furthermore, if $\xi_E$ is a characteristic element of $X(E/F_\infty)$, we have
\[\ord_{T=0}\big(\Phi_{\reg_{F_n}}(\xi_E)\big) =\left\{
                                         \begin{array}{ll}
                                            p^n-1, & \mbox{in case (a) of Proposition \ref{DD CS},} \\
                                            p^n, & \mbox{in case (b) of Proposition \ref{DD CS}.}
                                         \end{array}
                                       \right
.\]
Here $F_n = \Q(\mu_{p^n}, \sqrt[p^n]{q})$.
\ep

\bpf
As seen in Proposition \ref{DD CS}, $X(E/F_\infty)$ is a free $\Zp\ps{\Gal(F_\infty/\Q(\mu_{p^\infty}))}$-module of rank 1, and so in particular, it belongs to $\M_H(G)$. By \cite[Corollary 4.3]{CFKS}, we have
\[\corank_{\Zp}\big(\Sel(E/F_n)\big) \geq \left\{
                                         \begin{array}{ll}
                                            p^n-1, & \mbox{in case (a) of Proposition \ref{DD CS},} \\
                                            p^n, & \mbox{in case (b) of Proposition \ref{DD CS}.}
                                         \end{array}
                                       \right.\]
On the other hand, writing $H_n=\Gal(F_\infty/\mF_n)$ for $\mF_n = \Q(\mu_{p^\infty}, \sqrt[p^n]{q})$, and taking Proposition \ref{DD CS} and Remark \ref{hypotheses cd} into account, we may apply \cite[Lemma 2.3]{HungLim} to conclude that
\[ \rank_{\Zp}\big(X(E/F_\infty)_{H_n}\big) = p^n.\]
The argument in \cite[Theorem 3.1, Lemma 3.4]{HV} yields a map
\[ X(E/F_\infty)_{H_n} \lra X(E/\mF_n)\]
which has finite cokernel and whose kernel is of $\Zp$-rank $1$ or $0$ accordingly to case (a) or case (b) of Proposition \ref{DD CS}. In particular, this forces
\[\rank_{\Zp}\big(X(E/\mF_n)\big)=\corank_{\Zp}\big(\Sel(E/F_n)\big) = \left\{
                                         \begin{array}{ll}
                                            p^n-1, & \mbox{in case (a) of Proposition \ref{DD CS},} \\
                                            p^n, & \mbox{in case (b) of Proposition \ref{DD CS}.}
                                         \end{array}
                                       \right.\]
It then follows from this and Lemma \ref{semisimple lemma} that $X(E/\mF_n)$, as a $\Zp\ps{\Gal(\mF_n/F_n)}$-module, satisfies Conjecture \ref{semisimple conj}. The required conclusion of the proposition now follows from an application of Theorem \ref{main: pro-p}.
\epf

Finally, we relate the order of vanishing of the characteristic element to the order of zero of the Hasse-Weil $L$-functions $L(E/F_n,s)$ of $E$ at $s=1$. In according to the non-commutative Iwasawa main conjecture and the BSD-conjecture, one expects an equality between the two quantities. We can at least establish an inequality here.

\bp \label{CS DD}
Retain settings in Proposition \ref{DD CS}.
 If $\xi_E$ is a characteristic element of $X(E/F_\infty)$, then
 \[\ord_{T=0}\big(\Phi_{\reg_{F_n}}(\xi_E)\big)  \leq \ord_{s=1} \big(L(E/F_n,s)\big).\]
 \ep

\bpf
Suppose that we are in case (a) of Proposition \ref{DD CS}.
By \cite[Proposition A.38]{DD}, we have
\[ \ord_{s=1} \big(L(E/F_n,s)\big)\geq p^n-1.\]
Combining this with Proposition \ref{uncond}, we obtain the required inequality. Case (b) of Proposition \ref{DD CS} can be proven similarly by appealing to \cite[Proposition A.41]{DD} in place of \cite[Proposition A.38]{DD}.
\epf

Some examples of elliptic curves and primes $(p,q)$, where the preceding propositions can be applied, are $E = 11A3$ with $p= 3, q= 11$ and $E=38B1$ with $p= 3, q= 2$. For more examples of such elliptic curves and pairs of primes $(p,q)$, we refer readers to \cite[pp. 252-253]{DD}.

\medskip
In case (a) of Proposition \ref{DD CS}, by appealing to the work of Darmon-Tian \cite{DT}, the inequality in the preceding proposition can be improved to an equality under certain extra assumptions.

\bc \label{DarmonTian}
Let $E$ be an elliptic curve defined over $\Q$ with good ordinary reduction at $p$, and let $q$ be a multiplicative reduction prime of $E$. Suppose that all of the following statements holds.
\begin{enumerate}
  \item[$(a)$] The prime $q$ is inert in $\Q(\mu_p)/\Q$ and $\Sel(E/\Q(\mu_{p^\infty})) = 0$.
  \item[$(b)$] $\Gal(\Q(E[p]/\Q)\cong\mathrm{GL}_2(\Z/p\Z)$.
  \item[$(c)$] \cite[Conjecture 1.7]{DT} is valid.
\end{enumerate}
If $\xi_E$ is a characteristic element for $X(E/F_\infty)$, then
 \[\ord_{T=0}\big(\Phi_{\reg_{F_n}}(\xi_E)\big)  =\ord_{s=1} \big(L(E/F_n,s)\big).\]
 Furthermore, if $\rho_n$ is the Artin representation defined as in Subsection \ref{false-tate subsec}, we have
 \[\ord_{T=0}\big(\Phi_{\rho_n}(\xi_E)\big)  =\ord_{s=1} \big(L(E/\Q,\rho_n, s)\big). \]
 \ec

\bpf
Under the hypotheses of the corollary, Darmon and Tian have showed that $\ord_{s=1} \big(L(E/F_n,s)\big) = p^n-1$ (cf.\ \cite[Theorem 1.9]{DT}), and so the equality follows from combining this latter observation with Proposition \ref{uncond}. For the second equality, we first note that $\ord_{T=0}\big(\Phi_{\rho_n}(\xi_E)\big)$ is a non-negative integer by Theorem \ref{ArtinFT}. Moreover, since $X(E/F_\infty)$ is a free $\Zp\ps{\Gal(F_\infty/\Q(\mu_{p^\infty}))}$-module of rank 1, it follows from a parity result of Coates-Fukaya-Kato-Sujatha \cite[Theorem 4.6]{CFKS} that $\ord_{T=0} ~\big(\Phi_{\rho_n}(\xi_E)\big)=s_{E,\rho_n}$ is an odd integer. By (\ref{Phi rank}), we have
\begin{equation}\label{twist eqn}
  p^{n-1}(p-1)\ord_{T=0} ~\big(\Phi_{\rho_n}(\xi_E)\big) =  \corank_{\Zp} \big(\Sel(E/F_n)\big) - \corank_{\Zp} \big(\Sel(E/L_n)\big),
\end{equation}
   where $L_n =\Q(\mu_{p^{n}}, q^{p^{n-1}})$. Since $F_{n-1}\subseteq L_n$, one plainly has
   \[ \corank_{\Zp} \big(\Sel(E/L_n)\big) \geq \corank_{\Zp} \big(\Sel(E/F_{n-1})\big) = p^{n-1}-1. \]
   Putting this into (\ref{twist eqn}), we obtain
   \[  p^{n-1}(p-1)\ord_{T=0} ~\big(\Phi_{\rho_n}(\xi_E)\big) \leq p^n-1-(p^{n-1}-1) = p^{n-1}(p-1)  \]
   which in turn implies that
   \[ \ord_{T=0} ~\big(\Phi_{\rho_n}(\xi_E)\big)\leq 1.\]
   Since we have seen above that this integer is non-negative and odd, we may conclude that
   \[ \ord_{T=0} ~\big(\Phi_{\rho_n}(\xi_E)\big) = 1.\]
   On the analytic side, a combination of \cite[Theorems 1.8 and 1.9]{DT} and \cite[Theorem A.38]{DD} yields
   \[ \ord_{s=1} \big(L(E/\Q,\rho_n, s)\big) =1. \]
   This thus establishes the final equality and completes the proof of the corollary.
\epf

Finally, we note that Darmon-Tian has given many examples (see \cite[Section 3]{DT}), where hypotheses (a) and (b) are satisfied. Therefore, Corollary \ref{DarmonTian} applies to these elliptic curves (modulo \cite[Conjecture 1.7]{DT}).

\subsection{Multiplicative consideration}

We consider the analogue situation for an elliptic curve with multiplicative reduction. We first consider the split multiplication reduction situation.

\bl \label{Lee}
Let $E$ be an elliptic curve defined over $\Q$ with split multiplicative reduction at $p$. Let $q$ be either $p$ or a prime such that $E$ has non-split multiplicative reduction at the primes of $\Q(\mu_p)$ above $q$. Suppose that $\Sel(E/\Q(\mu_{p^\infty}))$ is finite.

Then $X(E/F_\infty)$ is a finitely generated $\Zp\ps{\Gal(F_\infty/\Q(\mu_{p^\infty}))}$-module of rank 1, where $F_\infty = \Q(\mu_{p^{\infty}},\sqrt[p^\infty]{q})$. \el

\bpf
This can be proven via a similar argument to that in \cite[Theorem 3.1]{HV}.
\epf

To continue, we recall a terminology of Lee \cite[Definition 1.10]{Lee}. For a given $p$-powerfree integer $m>1$, the pair $(p,m)$ is said to be amenable if either $p|m$ or $p | m^{p-1}- 1$ but $p^2\nmid m^{p-1}- 1$.

\bp \label{Lee2}
Retain settings in Lemma \ref{Lee}. Set $F_n= \Q(\mu_{p^n},\sqrt[p^n]{q})$. Let $\xi_E$ is a characteristic element of $X(E/F_\infty)$. If the pair $(p,q)$ is amenable, then
 \[\ord_{T=0}\big(\Phi_{\reg_{F_n}}(\xi_E)\big)  \leq \ord_{s=1} \big(L(E/F_n,s)\big)\]
 for $n\geq 1$.
  \ep

\bpf
 By \cite[Theorem 1.11]{Lee} and noting Lemma \ref{Lee}, we have
\begin{equation}\label{AMSel}
\corank_{\Zp}\big(\Sel(E/F_n)\big) = p^n-1
\end{equation}
for all $n\geq 1$. Building on this observation, we may apply a similar argument to that in Proposition \ref{uncond} to deduce that $X(E/\mF_n)$, as a $\Zp\ps{\Gal(\mF_n/F_n)}$-module, satisfies Conjecture \ref{semisimple conj}, where $\mF_n = \Q(\mu_{p^\infty}, \sqrt[p^n]{q})$. By an application of
Theorem \ref{main: pro-p}, we have
\begin{equation}\label{MultPhi}
\ord_{T=0}\big(\Phi_{\reg_{F_n}}(\xi_E)\big)  = \corank_{\Zp}\big(\Sel(E/F_n)\big).
\end{equation}
On the other hand, by a parity result of Lee \cite[Theorem 1.8 and Corollary 6.5]{Lee}, we have
\begin{equation}\label{AM1}
  w(E,\rho_n) = -1
\end{equation}
for all $n\geq 1$, where $\rho_n$ is the Artin representation defined as in Subsection \ref{false-tate subsec} and $w(E,\rho_n)$ is the root number in the sense of \cite{CFKS, Lee}. Artin formalism of the Hasse-Weil $L$-functions yields
\[ L(E/F_n,s) = \prod_{\psi} L(E/\Q, \psi, s)^{n_\psi},\]
where $\psi$ runs through all irreducible representations of $\Gal(F_n/\Q)$ and $n_\psi$ is the dimension of $\psi$. By (\ref{AM1}), we have
$\ord_{s=1}\big(L(E/\Q, \rho_n,s)\big)\geq 1$. Since $\rho_n$ has dimension $p^{n-1}(p-1)= p^n-p^{n-1}$ and $L(E/\Q, \psi, s)$ is holomorphic (cf. \cite[Theorem 14]{Do}), we have
\[ \ord_{s=1}\big(L(E/F_n,s)\big)\geq p-1 + p^2-p+\cdots + p^n-p^{n-1} = p^n-1. \]
Combining this with (\ref{AMSel}) and (\ref{MultPhi}), we obtain the conclusion of the proposition.
\epf

Finally, when the elliptic curve $E$ has non-split multiplicative reduction at $p$, we have the following analogous result.

\bp \label{LeeNonsplit}
Let $E$ be an elliptic curve defined over $\Q$ with non-split multiplicative reduction at $p$, and let $q$ be a multiplicative reduction prime of $E$. Suppose that one of the following statements holds.
\begin{enumerate}
  \item[$(a)$] $E$ has split multiplicative reduction at the prime $q$ with $q$ being inert in $\Q(\mu_p)/\Q$ and $\Sel(E/\Q(\mu_{p^\infty}))$ is finite.
  \item[$(b)$] $E$ has non-split multiplicative reduction at the primes of $\Q(\mu_p)$ above $q$ and $\Sel(E/\Q(\mu_{p^\infty})) = \Qp/\Zp$ with a trivial $\Ga$-action.
\end{enumerate}
Writing $F_n= \Q(\mu_{p^n},\sqrt[p^n]{q})$, we have
 \[\ord_{T=0}\big(\Phi_{\reg_{F_n}}(\xi_E)\big)  \leq \ord_{s=1} \big(L(E/F_n,s)\big)\]
 for $n\geq 1$.\ep

\bpf
This is proven similarly to that in Proposition \ref{Lee2}.
\epf

We refer readers to \cite[Section 7, Table 1]{Lee} for numerical examples, where Propositions \ref{Lee2} and \ref{LeeNonsplit} apply. Finally, we mention that since we do not have the analogous result of Darmon-Tian in the multiplicative reduction case, we are not able to establish equality as in Corollary \ref{DarmonTian}.

\end{document}